\documentclass{article}

\usepackage{geometry}
\usepackage{amsmath,amssymb}
\usepackage{changepage}
\usepackage[utf8x]{inputenc}
\usepackage{textcomp,marvosym}
\usepackage{cite}
\usepackage{nameref,hyperref}
\usepackage{array}
\usepackage{algorithm}
\usepackage[noend]{algpseudocode}
\usepackage{multicol}

\makeatletter
\renewcommand{\@biblabel}[1]{\quad#1.}
\makeatother

\usepackage{tikz}
\usetikzlibrary{decorations.markings}
\usetikzlibrary{decorations.pathreplacing}
\usetikzlibrary{spy}
\usetikzlibrary{positioning}
\usetikzlibrary{backgrounds}
\usetikzlibrary{hobby,decorations}
\usetikzlibrary{shapes,snakes}
\usetikzlibrary{shapes.geometric}
\usetikzlibrary{shadows}
\usetikzlibrary{intersections}
\usetikzlibrary{matrix}
\tikzset{   invisible/.style={opacity=0,text opacity=0},     visible on/.style={alt={#1{}{invisible}}},     alt/.code args={<#1>#2#3}{%
     \alt<#1>{\pgfkeysalso{#2}}{\pgfkeysalso{#3}} },}
\tikzset{cross/.style={cross out, draw=black, minimum size=2*(#1-\pgflinewidth), inner sep=0pt, outer sep=0pt},cross/.default={1pt}}
\tikzset{
  treenode/.style = {shape=rectangle, rounded corners,
                     draw, align=center,
                     top color=white, bottom color=blue!20},
  root/.style     = {treenode, font=\Large, bottom color=red!30},
  env/.style      = {treenode, font=\ttfamily\normalsize},
  dummy/.style    = {circle,draw},
  leaf/.style     = {circle, draw}
}

\usepackage{pgfplots}

\definecolor{blue-violet}{rgb}{0.54, 0.17, 0.89}
\usepackage{wrapfig}
\usepackage{color, colortbl}
\definecolor{Gray}{gray}{0.2}
\definecolor{aliceblue}{rgb}{0.94, 0.97, 1.0}
\definecolor{airforceblue}{rgb}{0.36, 0.54, 0.66}
\definecolor{antiquebrass}{rgb}{0.8, 0.58, 0.46}
\definecolor{cambridgeblue}{rgb}{0.64, 0.76, 0.68}
\definecolor{blue-violet}{rgb}{0.54, 0.17, 0.89}
\definecolor{blush}{rgb}{0.87, 0.36, 0.51}
\definecolor{carmine}{rgb}{0.59, 0.0, 0.09}
\definecolor{burgundy}{rgb}{0.5, 0.0, 0.13}
\definecolor{asparagus}{rgb}{0.53, 0.66, 0.42}
\definecolor{beaublue}{rgb}{0.74, 0.83, 0.9}

\usepackage{dsfont}

\newcommand{\1}{\mathds{1}}

\usepackage[sort=def,hyperfirst=false]{glossaries}
\newglossary*{disease}{\normalsize{Parameters of the disease model}}
\newglossary*{pdmp}{\normalsize{Parameters of the controlled PDMP/POMDP model}}
\newglossary*{pomcp}{\normalsize{Parameters of the POMCP algorithm}}
\newglossary*{acronyms}{\normalsize{Acronyms}}

\makeglossaries 

\newglossaryentry{glambda}{
type=disease,
name=\ensuremath{\lambda},
description={Disease risk function}
}

\newglossaryentry{gPhi}{
type=disease,
name=\ensuremath{\Phi},
description={Aggressiveness of the disease/treatment efficiency}
}

\newglossaryentry{gQ}{
type=disease,
name=\ensuremath{Q},
description={Transition kernel that dictates how the patient's state evolves at relapses}
}

\newglossaryentry{gm}{
type=disease,
name=\ensuremath{m},
description={The mode $m$ corresponds to the overall condition of the patient ($m=0$: remission, $m=1$: disease 1, $m=2$: disease 2, $m=3$: death of the patient)}
}

\newglossaryentry{gzeta}{
type=disease,
name=\ensuremath{\zeta},
description={The level of the marker}
}

\newglossaryentry{gzeta0}{
type=disease,
name=\ensuremath{\zeta_0},
description={The nominal marker value}
}

\newglossaryentry{gD}{
type=disease,
name=\ensuremath{D},
description={The death marker level}
}

\newglossaryentry{gu}{
type=disease,
name=\ensuremath{u},
description={The duration since when the patient has been in the current condition}
}

\newglossaryentry{gy}{
type=disease,
name=\ensuremath{y = \zeta+{\epsilon}},
description={denotes a noisy observation of the marker level}
}

\newglossaryentry{gomega}{
type=disease,
name=\ensuremath{\omega = (y, t)},
description={denotes a complete observation, where $t$ is the time since the beginning of the patient follow-up}
}

\newglossaryentry{gl}{
type=disease,
name=\ensuremath{\ell},
description={The treatment currently applied: $\ell=a$: Treatment efficient for disease 1, $\ell=b$: Treatment efficient for disease 2 and $\ell=\emptyset$: No treatment}
}

\newglossaryentry{gr}{
type=disease,
name=\ensuremath{r},
description={Delay untill next visit}
}

\newglossaryentry{gcv}{
type=disease,
name=\ensuremath{C_V},
description={Visit cost}
}

\newglossaryentry{gkappa}{
type=disease,
name=\ensuremath{\kappa},
description={Non-negative scale factor}
}

\newglossaryentry{gbeta}{
type=disease,
name=\ensuremath{\beta},
description={Penalty for unnecessary treatment}
}

\newglossaryentry{gM}{
type=disease,
name=\ensuremath{M},
description={Death cost}
}

\newglossaryentry{gX}{
type=pdmp,
name=\ensuremath{\left(X_t\right)_{0\leq t \leq H}},
description={Controlled PDMP/POMDP process}
}

\newglossaryentry{gY}{
type=pdmp,
name=\ensuremath{\left(Y_n\right)_{n\geq 0}},
description={PDMP/POMDP observation sequence}
}

\newglossaryentry{gs}{
type=pdmp,
name=\ensuremath{s=(m, \zeta,u)},
description={denotes the complete state of the patient in the PDMP/POMDP model of the disease control problem}
}

\newglossaryentry{gd}{
type=pdmp,
name=\ensuremath{d=(\ell,r)},
description={A decision in the PDMP/POMCP model of the disease control model}
}

\newglossaryentry{gpi}{
type=pdmp,
name=\ensuremath{\pi},
description={Treatment policy generating a sequence $(d_n)$ of decisions where $d_n = \pi\left(Y_0,\ldots,Y_n\right)$}
}

\newglossaryentry{ghn}{
type=pdmp,
name=\ensuremath{h_n=\langle\omega_0 d_0 \omega_1\ldots d_{n-1} \omega_n \rangle},
description={denotes the sequence of past observations and decisions, from which one can construct (theoretically) the current {\em belief state}}
}

\newglossaryentry{gbs}{
type=pdmp,
name=belief or filter space,
description={The true state of a POMDP is not directly observed. However, exploiting past decisions and observations (gathered in a {\em history}, one can maintain a {\em belief state}, which is a probability distribution over possible states of the POMDP. The belief space/filter space is the set of all the possible belief states}
}

\newglossaryentry{gV}{
type=pdmp,
name=\ensuremath{V},
description={Value of the optimal policy $\pi^*$, minimizing the expected cost}
}

\newglossaryentry{theta}{
type=pomcp,
name=\ensuremath{\Theta},
description={Generic filter used in the POMCP adaptation to controlled PDMPs, to generate particles sets}
}

\newglossaryentry{thetap}{
type=pomcp,
name=\ensuremath{\Theta^p},
description={Particle filter used in the POMCP adaptation to controlled PDMPs, to generate particles sets}
}

\newglossaryentry{thetac}{
type=pomcp,
name=\ensuremath{\Theta^c},
description={Conditional filters used in the POMCP adaptation to controlled PDMPs, to generate particles sets}
}

\newglossaryentry{B}{
type=pomcp,
name=\ensuremath{B},
description={Set of particles in the generic filter}
}

\newglossaryentry{Bp}{
type=pomcp,
name=\ensuremath{B^p},
description={Sets of particles in the particle filter}
}

\newglossaryentry{Bc}{
type=pomcp,
name=\ensuremath{B^c},
description={Sets of particles in the conditional filter}
}

\newglossaryentry{grollout}{
type=pomcp,
name=rollout strategy,
description={Arbitrary strategy which is applied in POMCP from leaves of the exploration tree to simulate trajectories. These trajectories are used to obtain preliminary estimates of the value function at the leaves}
}

\newglossaryentry{galphap}{
type=pomcp,
name=\ensuremath{\alpha'},
description={Normalised version of the exploration parameter $\alpha$ of POMCMP from Algorithm \ref{algo:POMCP}}
}

\newglossaryentry{gnsearch}{
type=pomcp,
name=\ensuremath{n_{\text{search}}},
description={Number of loops executed in the main loop of the adapted POMCMP. Plays the role of the {\sc Timeout} function in Algorithm \ref{algo:POMCP}}
}

\newglossaryentry{gK}{
type=pomcp,
name=\ensuremath{K},
description={Minimal size of the support of the particle filter or number of particles sampled from conditional filter}
}

\newglossaryentry{gPrec}{
type=pomcp,
name=\ensuremath{\mathcal D},
description={Discretization step of the observation space in the POMCP algorithm adapted to the disease control problem (with continuous observation space)}
}

\setacronymstyle{long-short}
\newacronym{dp}{DP}{Dynamic Programming}
\newacronym{ifm}{IFM}{Intergroupe Francophone du My\'elome}
\newacronym{fomdp}{MDP}{fully observed Markov Decision Process}
\newacronym{mm}{MM}{Multiple Myeloma}
\newacronym{pdmp}{PDMP}{Piecewise Deterministic Markov Process}
\newacronym{pfs}{PFS}{Progression-Free Survival}
\newacronym{pomdp}{POMDP}{Partially Observed Markov Decision Process}
\newacronym{pomcp}{POMCP}{Partially Observed Monte-Carlo Planning}
\newacronym{rl}{RL}{Reinforcement Learning}
\begin{document}

\title{Medical follow-up optimization: A Monte-Carlo planning strategy\thanks{We acknowledge the support of the French Agence Nationale de la Recherche (ANR), under grant ANR-21-CE40-005 (project HSMM-INCA), of European Union’s Horizon 2020 research and innovation program (Marie Sklodowska-Curie grant agreement No 890462) and the support of MESO@LR-Platform at the University of Montpellier.}}     

\author{Benoîte de Saporta\footnote{IMAG, Univ Montpellier, CNRS, Montpellier, France}
\and Aymar Thierry d'Argenlieu\footnote{IMAG, Univ Montpellier, CNRS, Montpellier, France and IP Paris, Palaiseau, France}
\and Régis Sabbadin\footnote{Univ Toulouse, INRAE-MIAT, Toulouse, France}
\and Alice Cleynen\footnote{IMAG, Univ Montpellier, CNRS, Montpellier, France and John Curtin School of Medical Research, Australian National University, Canberra, ACT, Australia}}
\date{}

\maketitle      

\begin{abstract}
Designing patient-specific follow-up strategy is a crucial step towards personalized medicine in cancer. 
Tools to help doctors deciding on treatment allocation together with next visit date, based on patient preferences and medical observations, would be particularly beneficial.
Such tools should be based on realistic models of disease progress under the impact of medical treatments, involve the design of (multi-)objective functions that a treatment strategy should optimize along the patient's medical journey, and include efficient resolution algorithms to optimize personalized follow-up by taking the patient's history and preferences into account.
We propose to model cancer evolution with a Piecewise Deterministic Markov Process where patients alternate between remission and relapse phases with disease-specific tumor evolution. This model is controlled via the online optimization of a long-term cost function accounting for treatment side-effects, hospital visits burden and disease impact on the quality of life. Optimization is based on noisy measurements of blood markers at visit dates.
This optimization problem is extremely difficult. It has recently been modeled as an infinite dimensional continuous space Markov Decision Process, approximated by a discrete-space problem in order to be solved exactly.
Here, instead, we leverage the Partially-Observed Monte-Carlo Planning algorithm to solve the full continuous-time, continuous-state problem, taking advantage of the nearly-deterministic nature of cancer evolution. We show that this approximate solution approach of the exact model performs better than the counterpart exact resolution of the discrete model, while allowing for more versatility in the cost function model, hence a patient-specific follow-up.
Our findings in terms of modeling and our efficient simulation-based optimisation approach to produce follow-up strategies can efficiently and easily be adapted to a large number of other diseases, thus being useful to doctors and patients.
\end{abstract}

\newpage
\tableofcontents
\newpage

\section{Introduction}
In long-term diseases such as cancer, patients alternate between remission
and relapse phases and are monitored along time through non-invasive check-ups such as blood
samples \cite{sozzi2001analysis,hundt2007blood}. Based on these noisy indirect disease measurements of some markers, practitioners must decide on treatment allocation, sometimes with little knowledge on the process dynamics (e.g. aggressiveness of the relapse)
which may differ between patients \cite{schrohl2003tumor,sharma2009tumor,nagpal2016tumor}. Long retrospect of medical practice has allowed the definition of milestones to help practitioners in making follow-up decisions, but automated personalized criteria are yet to be defined to improve individual patient follow-up.

The ability to monitor patients in the least invasive manner, according to their personal preferences (more check-ups to enforce relapse detection, less aggressive treatments for better quality of life, etc) is a crucial step towards better care, but requires fine knowledge of diseases dynamics and reliable prediction algorithms. One of the main requirements for such task is the definition of a universal model adapted to patient-specific parameters that could describe in an exhaustive manner the possible consequences of the practitioner's decisions. Mathematical models have been developed to link the tumor markers to tumor sizes \cite{lutz2008cancer}, or to predict evolution of tumor growth from initial measurements \cite{xu2016mathematical,nicolo2020machine}, but online adaptive models predicting relapses and automating treatment strategies are still lacking. In particular, such a model should be able to reconcile the continuous time evolution of the disease, continuous values for the markers 
leading to any possible values within a given range, and the noisy observations at discrete visit dates. {A good candidate is the class of \gls{pdmp} \cite{Davis84, C17, R17}. Indeed, PDMPs are non diffusive hybrid stochastic processes that can handle both continuous and discrete variables and their interactions in continuous time. The only source of stochasticity comes from the jumps of the process. They are thus simple to simulate and easy to interpret. Controlled PDMPs allow continuous time dynamics on continuous (or hybrid discrete and continuous) state spaces with decisions taken in continuous time \cite{Davis93,dSDZ15}.}

\begin{figure}[t]
\centering
\includegraphics[width=5cm]{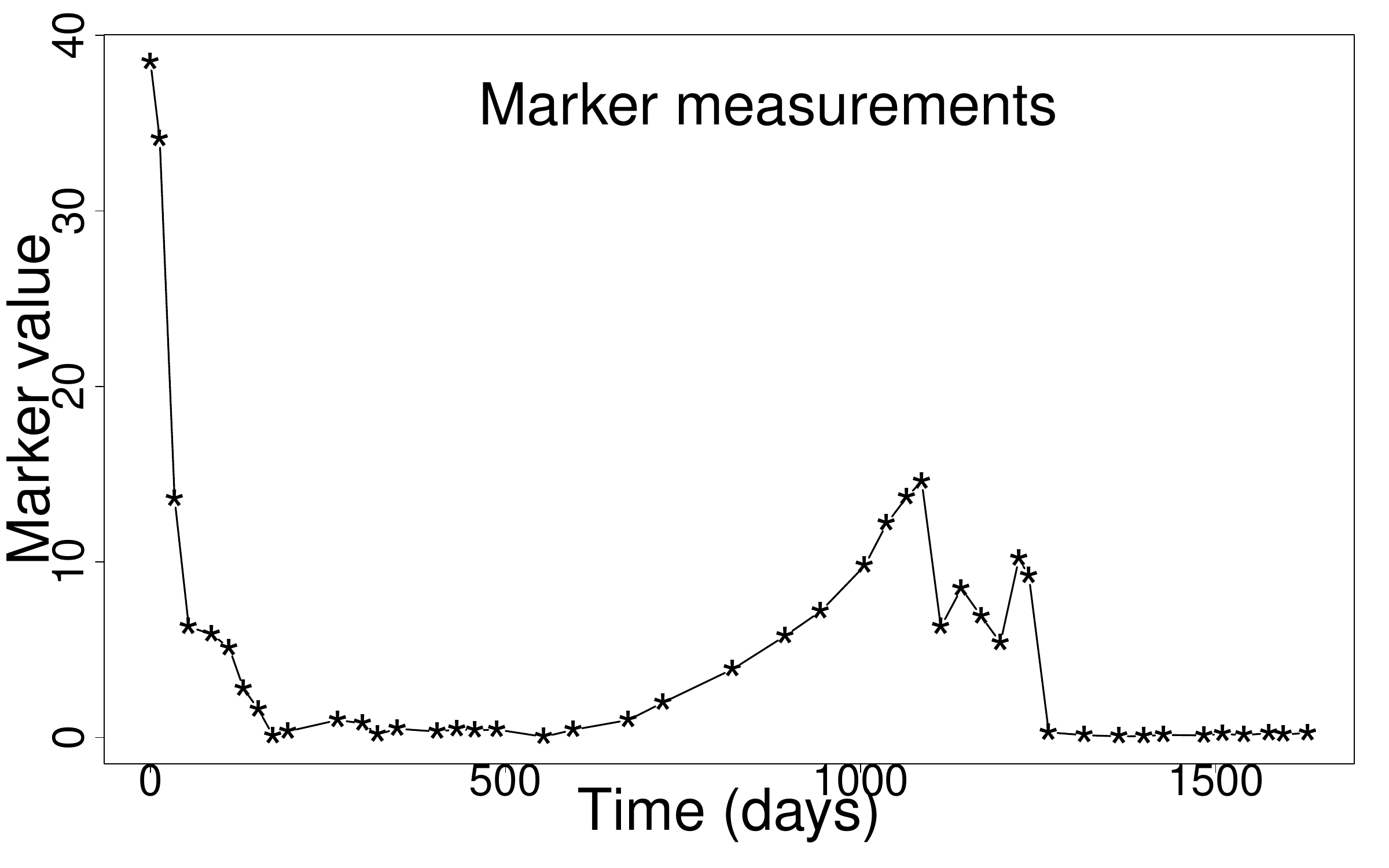} 
      \begin{tikzpicture}[,scale=0.72, every node/.style={scale=0.72}]
      \coordinate (t0) at (-2, 0);
      \coordinate(T1) at (-0.2, 0) ;
      \coordinate (d1) at (2.5,2);
      \coordinate (T2) at (4.5,0);
      \coordinate (d2) at (7,0);
      \coordinate (H) at (8,0);
      \draw[-] (t0)--(T1) 
      node[below=0.1,pos=1] {\small{relapse}};
      \draw (T1) to[out=0,in=-100] (d1)
      node[color=black, above=0.5,pos=0.9]{};
      \filldraw [gray] (T1) circle (2pt);
      \draw (d1) to[out=-90,in=-180] (T2) ;
      \draw[-] (T2)--(d2) node[below=0.1,pos=-0,color=black] {\small{remission}};
      \filldraw [gray] (T2) circle (2pt);
   \draw[dashed] (d2)--(H) node[below=0.1,pos=1,color=black] {$H$};
      \draw[-,color=black!50](-2,-0.1) -- (-2,0.1);
      \draw[-,color=black!50](-1,-0.1) -- (-1,0.1);
      \draw[-,color=black!50](0,-0.1) -- (0,0.1);
      \draw[-,color=black!50](1,0.15) -- (1,0.35);
      \draw[-,color=black!50](1.5,0.4) -- (1.5,0.6);
      \draw[-,color=black!50](2.5,1.9) -- (2.5,2.1)
      node[color=black, above=0.3]{\small{treatment}};
      \draw[->,color=black!50](2.5,2.4) -- (2.5,2.15);
      \draw[-,color=black!50](3,0.55) -- (3,0.75);
      \draw[-,color=black!50](3.5,0.15) -- (3.5,0.35);
      \draw[-,color=black!50](4.5,-0.1) -- (4.5,0.1);
      \draw[-,color=black!50](5,-0.1) -- (5,0.1);
      \draw[-,color=black!50](5.5,-0.1) -- (5.5,0.1);
      \draw[-,color=black!50](6,-0.1) -- (6,0.1);
      \draw[-,color=black!50](7,-0.1) -- (7,0.1);
      \draw[-,color=black!50](8,-0.1) -- (8,0.1);
      \draw[dashed] (-2,3)--(8,3) node[above=0.15cm,pos=1] {$D$};
      \draw[decorate, decoration = {brace}] (-2,0.2) -- (-0.2, 0.2)
      node[above=0.1,pos=0.5] {\small{risk $\lambda$}};
      \node at (1.4,1.5) {\small{disease}};
      \node at (0.9,1.2) {\small{aggressiveness}};
      \node at (1,0.9) {\small{$v$}};      
      \node at (3.5,1.5) {\small{treatment}};
      \node at (3.5,1.2) {\small{efficiency}};
      \node at (3.5,0.9) {\small{$v'$}};
    \end{tikzpicture}

\caption{\label{fig:patdata} {\bf Example of patient follow-up data, PDMP model. a)} Marker values are measured at each patient visits over a certain period of time. Data from the \textit{Intergroupe Francophone du Myélome} 2009 clinical trial, courtesy of the \textit{Centre de Recherche en Cancérologie de Toulouse}. {\bf b)} PDMP model, representation of the marker level of a patient. The risk function $\lambda$ controls the time to relapse, while parameters $v$ and $v'$ control the aggressiveness of the disease and the efficiency of the treatment respectively.}
\end{figure}
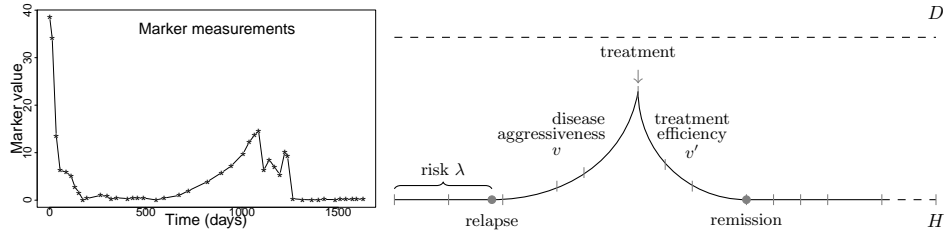

This paper is based on the analysis of a large cohort of \gls{mm} patient data from the \gls{ifm} 2009 clinical trial \cite{attal2017lenalidomide}. We assume that there is a single cancer marker that remains at a nominal threshold $\zeta_0$  throughout any remission phase, and that at patient relapse its level increases exponentially with multiple possible behaviors until treatment is administered, or a threshold $D$ is reached and the patient dies. To set up the context, we will assume that the study begins at time $t_0=0$ when the patient enters her first remission state, and we will denote $t_1$, $t_2,\ldots$, her visit dates, defined over time by the practitioner, until some time horizon $H$ is reached, or the patient dies. The time lapse between visits may not be constant, so that different patients may have different visit numbers and dates. More precisely, we will assume that at each visit time, the practitioner may choose to schedule the next visit in either $15$, $30$ or $60$ days. 
Such decision may be based on the previous and current marker measurements, which we denote $Y_0$, $Y_1, \dots$. Note that the exact value of the marker is hidden as measurements are corrupted by noise, and measurements are only collected at visit dates. Together with the next visit date, the practitioner may chose to modify the patient's current treatment, fixing it to one of the two available treatments, $a$ and $b$, or to no treatment at all, denoted $\emptyset$. An example of patient follow-up data is presented in Figure \ref{fig:patdata} a).

We model the underlying continuous-time dynamics of the patient health by a controlled 
PDMP \gls{gX} 
\cite{Davis93,DY95,Al01,CD13,dSDZ15,RTT17}. 
We propose to optimally control the process, that is to choose online the next treatment and visit date, based on present and past observations and decisions 
by minimising a cost function which is calibrated to balance the burden of deteriorated quality of life under treatment (including hospital visits) with the risk of dying from the disease. 

Previous work has focused on 
discretizing this problem in order to solve it approximately through \gls{dp} iterations \cite{CdS23}.  More specifically, the optimal control problem for the PDMP has first been expressed as a 
\gls{pomdp} \cite{BR11}. 
This step is simply done by considering decision dates as stages of the POMDP. Note that the time lapse between decisions is thus not constant, the continuous time dynamics is encoded in the specific parametrization of the transition kernel, and the POMDP still has a continuous state space, with continuous observation space.
The problem is then classically converted into a \gls{fomdp} \cite{BR11,CdS18}
on the {\em \gls{gbs}}. 
The filter process represents the probability distribution of the hidden values of the patient current state given the past and present observations. Second, the state space of the controlled PDMP has been discretized, so that an approximation of the filter process could be computed, charging only finitely many states. This approximate filter is called {\em conditional filter} in the sequel. 
Third, the belief space has been 
discretized in order to solve the MDP via dynamic programming iterations on a finite space.

In the current article, instead of discretizing the state and belief spaces, we use a {\em Monte-Carlo Tree Search} approach to (approximately) solve the controlled PDMP problem by simulation. 
More precisely, we propose an adaptation of the \gls{pomcp}
algorithm \cite{SV10}, originally designed to solve discrete time / finite state and observation spaces POMDP, to the case of controlled PDMP. The novel challenge is that controlled PDMP involve continuous time as well as continuous state and observation spaces.  
We show empirically that this simulation-based approach outperforms the discretization-based approach both in terms of computation time and quality of returned policies.
Thus, this approach is promising for providing an automated decision aiding tool for practitioners. 
 To our knowledge, this is the first method to address this challenging question. Current decision making is typically based on heuristic rules derived from expert clinical knowledge \cite{stewart2009treat,lonial2015treat,gerecke2016diagnosis}.

\section{Results}
\subsection{Controlled piecewise deterministic Markov processes form a universal class of versatile models for patient follow-up}
Despite the discrete-time acquisition of the marker measurements, we choose to model the dynamics of the patient's health by a continuous-time controlled Piecewise Deterministic Markov Process (PDMP). The formalism of PDMPs is both light and versatile \cite{Davis93,RTW12,C17,R17}. It allows to describe the dynamics of the disease 
with only three biologically relevant parameters: the disease's risk function, $\gls{glambda}$, that dictates how often the patient is likely to relapse, the aggressiveness of the disease $v$ that dictates how fast the marker level will increase during relapses, and the treatment efficiency $v'$ that dictates   how fast the marker level decreases under treatment. \\
 In the formalism of PDMP, this is formulated by an exponential flow $\gls{gPhi}$ which slope parameter ($v$ or $v'$) depends on patient condition and treatment, the risk function $\lambda$,
and a transition kernel $\gls{gQ}$, that dictates how the patient's state evolves at relapses, here preventing the marker values from jumping abruptly at patient condition changes. This is illustrated in Fig \ref{fig:patdata} b). 

We consider a common risk function (identical for all patients) which we allow to depend on the time since the last remission date as well as the cancer marker level. We assume that the aggressiveness of the disease can be patient-dependent. It is either high 
or low 
and we model this as two different diseases. This leads us to introduce two different treatments, each efficient for one of those diseases and slowing the progression of the other. It is also an option not to treat for a given period.

We introduce three variables $m,\zeta,u$, where the mode $\gls{gm}$ corresponds to the overall condition of the patient ($m=0$: remission, $m=1$: disease 1, $m=2$: disease 2, $m=3$: death of the patient), 
$\gls{gzeta}\in[\zeta_0,D]$ is the level of the marker, 
where $\gls{gzeta0}$ is the nominal value and $\gls{gD}$ the death level 
and $\gls{gu}\geq 0$ is the time since the last change of overall condition (added for technical reasons to deal with non-constant risk functions). 
The precise definition of the controlled PDMP and its parameters $(\lambda, \Phi, Q)$ are given in the Methods section. 
The complete state of the patient is thus encoded by $\gls{gs}$. We  denote $X_0, \dots, X_n$ the process values at the observation dates $t_0, \dots, t_n$.

The overall condition of the patient $m$, the level of the marker $\zeta$ and the relapse dates (together with the time $u$ since the last change of condition) are not directly observed and thus cannot be used by the clinician to select a treatment. At each visit of the patient to the medical center, we assume that the practitioner receives a noisy observation of the marker level $\gls{gy}$ where ${\epsilon}$ is some Gaussian noise. The practitioner also knows the time $t$ since the beginning of the patient follow-up. The complete observation available to the practitioner is thus encoded by $\gls{gomega}$. 
The practitioner also has access to an indicator that the patient is still alive as treatment and follow-up stop at the death of the patient.

Based on the collection of present and past measurements and decisions, the practitioner selects both a time delay $\gls{gr}$ until the next visit to the medical center and a treatment $\gls{gl}$ to hold until this next visit.
Note that in our framework measurements are only made at visit dates.
A decision is thus a pair $\gls{gd}$, where $\ell \in \{\emptyset,a,b\}$, and $r\in \{15,30,60\}.$ Given a fixed arbitrary decision policy, simulating controlled patient trajectories is easy: see Algorithm \ref{SimuPDMP} given in the Methods Section. 

\subsection{Cost functions encode the diverse impacts of treatment on the patient's quality of life}
For the practitioner, controlling the disease is equivalent to choosing the best available treatment as well as the best next visit date in order to minimize its impact on the patient's quality of life along time. %
Defining the impact of treatment on the quality of life is a difficult task as it will typically depend on the treatment's side effects, the number of visits, the burden of living with a disease and the remaining life expectancy.

This paper proposes a mathematical definition of the impact of the treatment on quality of life in terms of a cost function that takes into account those different aspects. For a decision $d=(\ell, r)$ comprising a treatment allocation $\ell$ and a time to next visit $r$, and for a current marker level $\zeta$ at time $t_k$, and future marker level $\zeta'$ at time $t_{k+1}=t_k+r$ , we define 
\begin{align}
c(\zeta, d, \zeta') = C_V  +\kappa | \zeta' -\zeta_0 |r +\beta r \1_{\{\zeta=\zeta_0,\ell\neq \emptyset\}}+ M  \1_{\{\zeta'=D\}},
\end{align}
where $\gls{gcv}$ is a visit cost, $\gls{gkappa}$ is non-negative scale factor penalizing high marker values, $\gls{gbeta}$ is a penalty for applying an unnecessary treatment and $\gls{gM}$ is the death cost. 

This cost function thus takes into account a visit cost, to prevent patients from undergoing too many screening tests, a cost depending on the marker value at the next visit, to encourage treatment and calibrate visit dates, a cost for degradation of quality of life due to treatments, in particular if they are not appropriate, and a cost for dying.

Calibrating cost parameters $C_V$, $\kappa$, $\beta$, and $M$ is a very difficult task, which is
allowed to be patient-dependent (some patients
may even express a wish to be sedated rather than undergo
very long and painful treatments), and treatment
strategies are bound to be parameters-dependent.

When cast as a controlled PDMP with this cost function, the practitioner's problem is mathematically equivalent to solving a (continuous state space) Partially Observable Markov Decision Process (POMDP) \cite{KLC98,CdS23}, which expected value optimisation can be stated as
\begin{align}
V = \inf_{\pi \in \Pi} \mathbb E_{X_0}^\pi \left[\sum_{n=0}^{N_\pi-1} c(X_n,d_n,X_{n+1}) \right], \label{eq:VF}
\end{align}
where $\gls{gV}$ is called the \emph{optimal policy value} and represents the lowest possible expected total cost, $\Pi$ is the set of admissible policies (yielding decisions depending only on current and past observations), $N_\pi$ is the patient-specific  number of visits within the time-horizon of the study when using policy $\gls{gpi}$, $d_n=\pi(Y_0,t_0,\ldots,Y_n,t_n)$ is the decision $(\ell,r)$ taken at the $n$-th visit date $t_n$ according to policy $\pi$, and \gls{gY} represents the marker observation process for the controlled-PDMP/POMDP. Solving this problem amounts to computing (a good approximation of) the optimal policy value and identifying an admissible policy $\pi^*$ that reaches (a value close to) the minimum. 

\subsection{Adapted Partially Observed Monte-Carlo Planning is particulartly well suited for controlled PDMPs}
The Partially Observed Monte-Carlo Planning (POMCP) algorithm \cite{SV10} is an efficient simulation-based algorithm that has been designed for real-time planning in large finite state-space POMDPs. In this paper we show that even-though it has not been designed to handle continuous state and observation spaces, we can adapt it to solve {controlled} PDMPs, thanks to their efficient simulation property, without resorting to the computation of complex integrals for computing transition probabilities.

The objective of POMCP is to reduce the complexity of dynamic programming, which requires the construction of the entire decision tree (including the probabilities of every possible outcome with every possible decision at every future time-point), by sampling the tree in a principled way so as to compute the current optimal action.
POMCP is thus an {\em online} algorithm, which re-estimates the optimal strategy at each new data acquisition\footnote{POMCP does not "forget" previously computed strategies, but updates them using new simulated samples after every new observation is received.}. 

The POMCP algorithm relies on two main properties.
The first one is the ability to simulate trajectories, so as to progressively build the decision tree and update filters \gls{theta} at every intermediate node $h$ of the tree. Recall that a filter is a probability distribution representing the (approximate) distribution of the current hidden state given the observations. The standard POMCP algorithm uses a specific family of simulation-based filters \gls {thetap} called \emph{particle filters} specified below. Filters are used to sample sets 
of plausible states.

The second property is the requirement to provide {\em estimates} of the expected value of the policy in leaves of the current exploration tree, in order to guide exploration and build the decision policy. 

In the Methods Section we detail the algorithm (Algorithm \ref{algo:POMCP}) and show that POMCP is particularly well suited for controlled PDMPs.  
We simply point out here why trajectories simulation and policy evaluation are particularly efficient in POMCP, in the case of a controlled PDMP.
\begin{enumerate}
    \item  Simulation is particularly straightforward with PDMPs  \cite{Gill77,dSDZ15,LTT20}, requiring only to simulate the jump times and exploit the deterministic behavior between jumps, see Algorithm \ref{SimuPDMP}. In our medical framework, it is made even more simple since only few jumps are allowed.
    When little knowledge is available about the underlying process, a classic approach is to resort to \textit{particle filters} $\Theta^p$ \cite{DM96}. A particle filter $\Theta^p$ at step $n$ is a discrete uniform probability distribution with finite support 
    \gls{Bp} (where $B^p$ may have repeated atoms). It is updated at step $n+1$ though simulations: states $s$ from  $B^p$  are updated through a one-step simulation to a new state $s'$, and selected to be added to $B^p$ if the simulated observation is close to the true one. 
    As an alternative filter to compare to, we propose to use a \textit{conditional filter} 
    \gls{thetac} derived from the exact filter (that is the conditional distribution of the hidden state given the observations) from \cite{CdS23}. The exact conditional filter is updated through a recurrence formula involving ratios of integrals over the state space. By discretizing the state space, one can construct the approximation $\Theta^c$ of the exact filter. Unlike the particle filter that has a dynamically changing support with a uniform mass function, this conditional filter has a fixed support (the discretized space) with changing mass functions that are updated through analytical ratios of weighted sums.
    \item To estimate the future expected cost at some node of the tree, POMCP requires to simulate many full trajectories from the current node to a leaf of the tree. This requires to apply an arbitrary strategy to pick actions at every future nodes for which a decision has not yet been optimized. This arbitrary strategy is called a \emph{\gls{grollout}} in the POMCP framework. The most naive rollout strategy consists in uniformly randomly selecting decisions from the decision set $\{\emptyset,a,b\}\times\{15,30,60\}$. We consider instead a mode-based rollout strategy, which consists in choosing action $\emptyset$ in mode $0$ (no treatment if the simulated patient is in remission), action $a$ in mode $1$ and $b$ in mode $2$ (most efficient treatment if the simulated patient has relapsed) and a fixed next visit date of $15$ days. This rollout strategy, while not being necessarily optimal (depending on the cost function it might be optimal not to treat at the beginning of a relapse, or to treat preventively when in remission), exploits knowledge of the cost function, hence yields better estimates of action costs at time $t$. Note also that this mode-based policy is not applicable for real patients, since their mode is not observed. It is only applicable to simulated patients.
 This is fine since POMCP's rollout strategy is only used through simulations to estimate costs.
\end{enumerate}

POMCP has a number of tuning parameters (number of simulations, number of particles in the filter, exploration vs exploitation rates) which are described (as well as the impact of varying their values) in the following section.

\subsection{Following up a patient with adapted POMCP is easy and fast in practice}
\begin{figure}[htp]
  \begin{center}
      \begin{tikzpicture}[,scale=1, every node/.style={scale=1},
    grow                    = right,
    sibling distance        = 6em,
    level distance          = 10em,
    edge from parent/.style = {draw, -latex},
    every node/.style       = {font=\footnotesize},
    sloped
  ]
      \node at (-2.75, 2.35){{\bf a)}};
      \draw[rounded corners]  (-2.5, -1) rectangle (9, 2.5) node (Blocka) {};
      \coordinate (t0) at (-2, 0);
      \coordinate (H) at (8,0);
   \draw[dashed,->] (t0)--(H) node[below=0.1,pos=1,color=black] {$H$};
      \filldraw [gray] (-2,0.9) circle (2pt)node[below=0.1]{$y_0$};
      \filldraw [gray] (-1,1.4) circle (2pt)node[below=0.1]{$y_1$};
      \filldraw [gray] (0,1) circle (2pt)node[below=0.1]{$y_2$};
      \filldraw [gray] (1.5,0.9) circle (2pt)node[below=0.1]{$y_3$};
      \filldraw [gray] (2.5,1.9) circle (2pt)node[below=0.1]{$y_4$};
      \filldraw [gray] (3,1.8) circle (2pt)node[below=0.1]{$\dots$};
      \filldraw [gray] (3.5,1.1) circle (2pt)node[below=0.1]{$\dots$};
      \filldraw [gray] (4.5,1.7) circle (2pt)node[right=0.1,text width=3.5cm,,color=black]{\small{new marker acquisition $\omega_n=(y_n,t_n)$}};
      \draw[-,color=black!50](-2,-0.1) -- (-2,0.1)node[below=0.1]{$t_0$};
      \draw[-,color=black!50](-1,-0.1) -- (-1,0.1)node[below=0.1]{$t_1$};
      \draw[-,color=black!50](0,-0.1) -- (0,0.1)node[below=0.1]{$t_2$};
      \draw[-,color=black!50](1.5,-0.1) -- (1.5,0.1)node[below=0.1]{$t_3$};
      \draw[-,color=black!50](2.5,-0.1) -- (2.5,0.1)node[below=0.1]{$t_4$};
      \draw[-,color=black!50](3,-0.1) -- (3,0.1);
      \draw[-,color=black!50](3.5,-0.1) -- (3.5,0.1)node[below=0.15]{$\dots$};
      \draw[-,color=black!50](4.5,-0.1) -- (4.5,0.1);
      \node at (-2.75, -2.05){{\bf b)}};
      \draw[rounded corners]  (-2.5, -6) rectangle (9, -1.9) node (Blockb) {};
    \draw[->]  (2.75,-1) -- (2.75,-1.9);
    \node at (3.1, -1.45) {$\omega_n$} ;
    \begin{axis}[ybar interval, ymax=55,ymin=0, width=0.25\linewidth, yticklabels={},at={(-0.015\linewidth,-0.325\linewidth)},xtick=\empty,hide axis]
    \addplot coordinates { (0, 5) (2.5,5) (5, 35) (7.5,40) (10, 50) (12.5,20) (15, 15) (17.5,15) (20, 13) (22.5, 0) };
    \end{axis}
    \begin{axis}[xbar interval, xmax=55,xmin=0, width=0.25\linewidth, yticklabels={},at={(-0.12\linewidth,-0.28\linewidth)},ytick=\empty,hide axis]
    \addplot coordinates { (10, 0) (30,2.5) (50, 5) (40,7.5) (15, 10) (5,12.5) (5, 15) (3,17.5) (0, 20) (0, 22.5) };
    \end{axis}
    \begin{axis}[ybar interval, ymax=55,ymin=0, width=0.25\linewidth, yticklabels={},at={(0.41\linewidth,-0.325\linewidth)},xtick=\empty,hide axis]
    \addplot coordinates { (0, 3) (2.5,5) (5, 8) (7.5,20) (10, 45) (12.5,50) (15, 20) (17.5,45) (20, 20) (22.5, 10) };
    \end{axis}
    \begin{axis}[xbar interval, xmax=55,xmin=0, width=0.25\linewidth, yticklabels={},at={(0.325\linewidth,-0.28\linewidth)},ytick=\empty,hide axis]
    \addplot coordinates { (30, 0) (50,2.5) (40, 5) (30,7.5) (10, 10) (5,12.5) (5, 15) (0,17.5) (10, 20) (25, 22.5) (0,25) };
    \end{axis}
    \draw[->] (2.8,-3.4)--(3.8,-3.4) node[above=0.1,pos=0.5,color=black] {update};
    \draw[->] (-1.8,-5)--(-1.8,-2.5);
        \node at (-1.95, -3.75) {\small{$\zeta$}};
    \draw[->] (-1.9,-4.9)--(1.9,-4.9) node[below=0.1,pos=0.5,color=black] {\small{$u$}};
    \draw[->] (4.85,-5)--(4.85,-2.5);
    \node at (4.7, -3.75) {\small{$\zeta$}};
    \draw[->] (4.75,-4.9)--(8.4,-4.9) node[below=0.1,pos=0.5,color=black] {\small{$u$}};
    \node at (0,-5.5) {$\Theta_{n-1} $};
    \node at (6.9,-5.5) {$\Theta_{n} $};
    \draw[->]  (2.75,-6) -- (2.75,-6.9);
    \node at (3.5, -6.5) {$(h_n, \Theta_n) $} ;
    \draw[->]  (-3.5,-6.9) -- (-3.5, 0.75) -- (-2.5,0.75);
    \node at (-3.75, -3.25) {$d_n$} ;
    \draw[rounded corners]  (-4.05, -13) rectangle (9, -6.9) node (Blockc) {};
  \node at (-4.25, -7.1){{\bf c)}};
   \node [root] at (-3, -11) {$h_n$\\ $\langle N, V\rangle$}
     [level distance =25mm]
     child { node [env] {$h_nd_1^n$\\ $\langle N, V\rangle$}
         [level distance =90mm] child { node [leaf] {$c$}
            edge from parent[dashed] node [above] {rollout} }
      edge from parent node [below] {$d_1^n$}
      }
    child { node [env] {$h_nd_2^n$\\ $\langle N, V\rangle$}
      [level distance =30mm]
      child { node [env] {$h_nd_2^n\omega_1^{n+1}$\\ $\langle N, V\rangle$}
        [level distance =60mm] child { node [leaf] {$c$}
            edge from parent[dashed] node [above] {rollout} }
        edge from parent node [below] {$\omega_1^{n+1}$} }
      child { node [env] {$h_nd_2^n\omega_2^{n+1}$\\ $\langle N, V\rangle$}
        child { node [env] {$h_nd_2^n\omega_1^{n+1}d_1^{n+1}$\\ $\langle N, V\rangle$}
            child { node [leaf] {$c$}
                edge from parent[dashed] node [above] {rollout} }
            edge from parent node [above] {$d_1^{n+1}$}
            }
        child { node [env] {$h_nd_2^n\omega_1^{n+1}d_2^{n+1}$\\ $\langle N, V\rangle$}
            child { node [leaf] {$c$}
                edge from parent[dashed] node [above] {rollout} }
            edge from parent node [above] {$d_2^{n+1}$} }
        edge from parent node [above] {$\omega_2^{n+1}$}
        }
    edge from parent node [above] {$d_2^n$}
              };
\end{tikzpicture}
\end{center}
\caption{\label{fig:procedure} {\bf Practice of patient follow-up. a)} At each new visit the patient has a new marker measurement, and the practitioner receives a new observation $\omega_n=(y_n,t_n)$. {\bf b)} The filter is updated with the new observation, either through particle rejection sampling (\textit{particle filter}) or via a recursion formula (\textit{conditional filter}). {\bf c)} The decision tree is partially explored via simulation through an adapted POMCP algorithm using the updated filter. The algorithm returns the optimal decision $d_n$, combination of a time to next visit (defining $t_{n+1}$) and treatment to allocate (influencing $y_{n+1}$).  }
\end{figure}
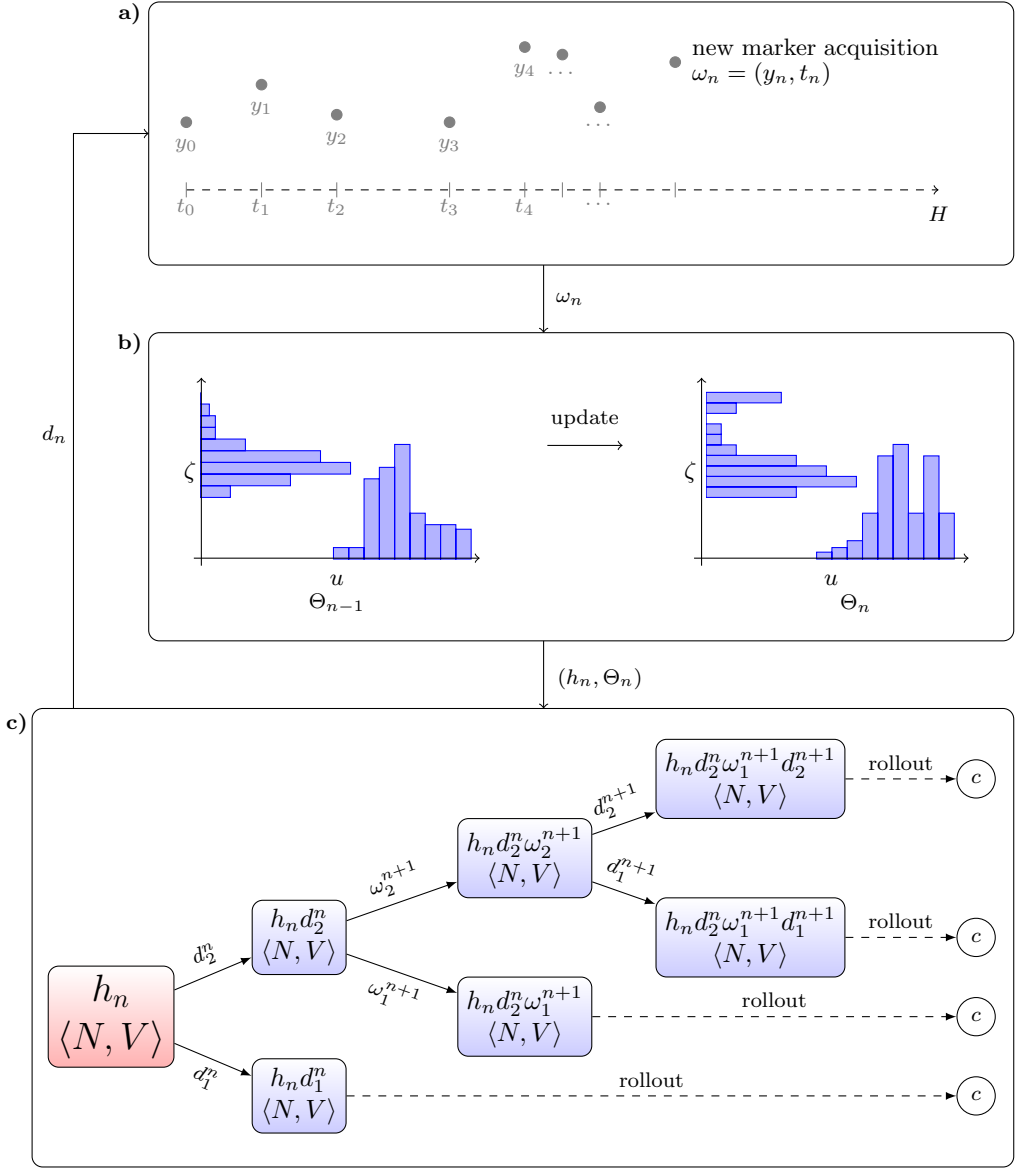

The previous paragraphs set the grounds for optimizing the long-term follow up of patients. In practice, we will assume a patient will enter the follow-up study once she enters the remission phase after an initial round of treatment. The practitioner may hence assume that her current state is known, \textit{i.e.} $s_0=(0,\zeta_0,0)$ and the initial value of both the particle and conditional filters is the Dirac mass at $s_0$. The initial observation is $\omega_0=(\zeta_0, t_0)$. The adapted POMCP algorithm is run to obtain the optimal decision $d_0$, which the practitioner can use (if she decides to) to allocate treatment and decide on the next visit date $t_1$. \\
At visit $n$, the patient will come back for some new marker measurement, so that the $n$-th observation value $\omega_n=(y_n,t_n)$ is obtained. The practitioner will have access to her full history, $h_n = \langle \omega_0 d_0 \omega_1 d_1 \cdots d_{n-1} \omega_{n} \rangle$ as well as her last belief filter, $\Theta^c_{n-1}$ or $\Theta^p_{n-1}$. An initial update of the filter is performed, either using the recursion formula for $\Theta^c_{n}$ from $\Theta^c_{n-1}$ and $\omega_n$, or by particle filtering through rejection sampling for $\Theta^p_{n}$ from $\Theta^p_{n-1}$ and $\omega_n$. The adapted POMCP algorithm is then ran to obtain the optimal current decision $d_n$, which the practitioner can use (or not) to allocate treatment and decide on the next visit date $t_{n+1}$. This is illustrated in Figure \ref{fig:procedure}.

\subsection{Adapted POMCP can be tuned to outperform dynamic programming}

The simulation study presented in this section has been conducted based on real data obtained from the \textit{Centre de Recherche en Cancérologie de Toulouse} (CRCT). Multiple myeloma (MM) is the second most common haematological malignancy in the world and is characterised by the accumulation of malignant plasma cells in the bone marrow. Classical treatments are based on chemotherapies, which, if appropriate, act fast and efficiently bring MM patients to remission in a few weeks. However almost all patients eventually relapse more than once and the five-year survival rate is around 50\%. 

We have obtained data from the \textit{Intergroupe Francophone du Myélome} 2009 clinical trial \cite{attal2017lenalidomide} which has followed 748 French MM patients from diagnosis to their first relapse on a standardized protocol for up to six years. At each visit a blood sample was obtained to evaluate the amount of monoclonal immunoglobulin protein in the blood, a marker for the disease progression. An example of patient dataset is given in Figure~\ref{fig:patdata}. 

Based on these data, we calibrated our PDMP model as described in the Methods section, and  we performed simulations to evaluate the performance of the POMCP strategy to select the combination of treatment and next visit date at each time point of the trajectories (these time-points being themselves selected by the algorithm). The performance of the approach was measured by a Monte-Carlo estimate of the expectation and confidence interval of its value as well as the runtime of the online computation of a complete trajectory. For each disease parameters' configuration $500$ simulations were performed to estimate these values. Codes and parameters are available at \url{https://github.com/acleynen/pomcp4pdmp} \cite{pomcp4pomdp}.\\

\subsubsection{Study 1: Impact of the parameters' values on POMCP's perfomance}
We evaluated the impact of the value of $6$ parameters on the performance: 
(i) the \textit{filter} chosen (conditional or particle), 
(ii) the \textit{rollout} procedure chosen, 
(iii) the exploitation/exploration tradeoff parameter \gls{galphap}, 
(iv) the number \gls{gnsearch} of simulations in the online POMCP procedure,
(v) the number \gls{gK} of initial states to sample from at each of the $n_{\text{search}}$ simulations, 
and 
(vi) the internal POMCP precision parameter \gls{gPrec} to select particles in the particle filter. 
Those parameters are described at length in the Methods Section, see Algorithm~\ref{algo:POMCP}. 
The github page \cite{pomcp4pomdp} contains tables with results for every sets of parameters' values that were tested. In this section we describe the most important results. \\
 We performed all parameter comparisons for both the conditional and the particle filters. The conditional particle filter is a discrete probability distribution on a finite fixed support $\gls{Bc}$ of size $184$ with $81$ states in condition $m=0$, $31$ states in condition $m=1$, $71$ states in condition $m=2$ and one state in condition $m=3$. The choice of these states is discussed in \cite{CdS23}. To adapt this filter to the POMCP environment, at each iteration $n$ we start by randomly sampling $K$ states $s$ from $B^c$ with the distribution given by $\Theta_n^c$. For the particle filter, this number $K$ directly corresponds to the number of particles in the filter, hence the size of $B^p$. Note that for the conditional filter the support $B^c$ does not change over time, whereas for the particle filter $B^p$ keeps the same size but possibly contains  different states at each iteration.

\paragraph{Mode-based rollout outperforms naive rollout.} We found that the uniform rollout procedure (selecting decisions randomly) produced very poor results compared to the mode-based rollout procedure 
and hence we only present results for the mode-based rollout policy here. 

\paragraph{POMCP is robust to the exploration/exploitation trade-off.} We found that the exploitation/exploration trade-off parameter had little influence on the overall performance, with the exception of extreme values ($\alpha'=0.99$, almost no exploration, and $\alpha'=0.2$, almost no exploitation, both leading to poorer performance). We also tried several adaptive strategies to select the value of $\alpha'$ depending on the confidence in the belief (measured in terms of entropy) which did not improve the results. The following results are therefore discussed for a fixed value of $\alpha'=0.5$, but the reader may refer to the Supplementary Information Tables \ref{table_resultats_1} and \ref{table_resultats_2} for additional results on these parameters. 

\begin{figure}[ht!]
\centering
\includegraphics[width=12.5cm]{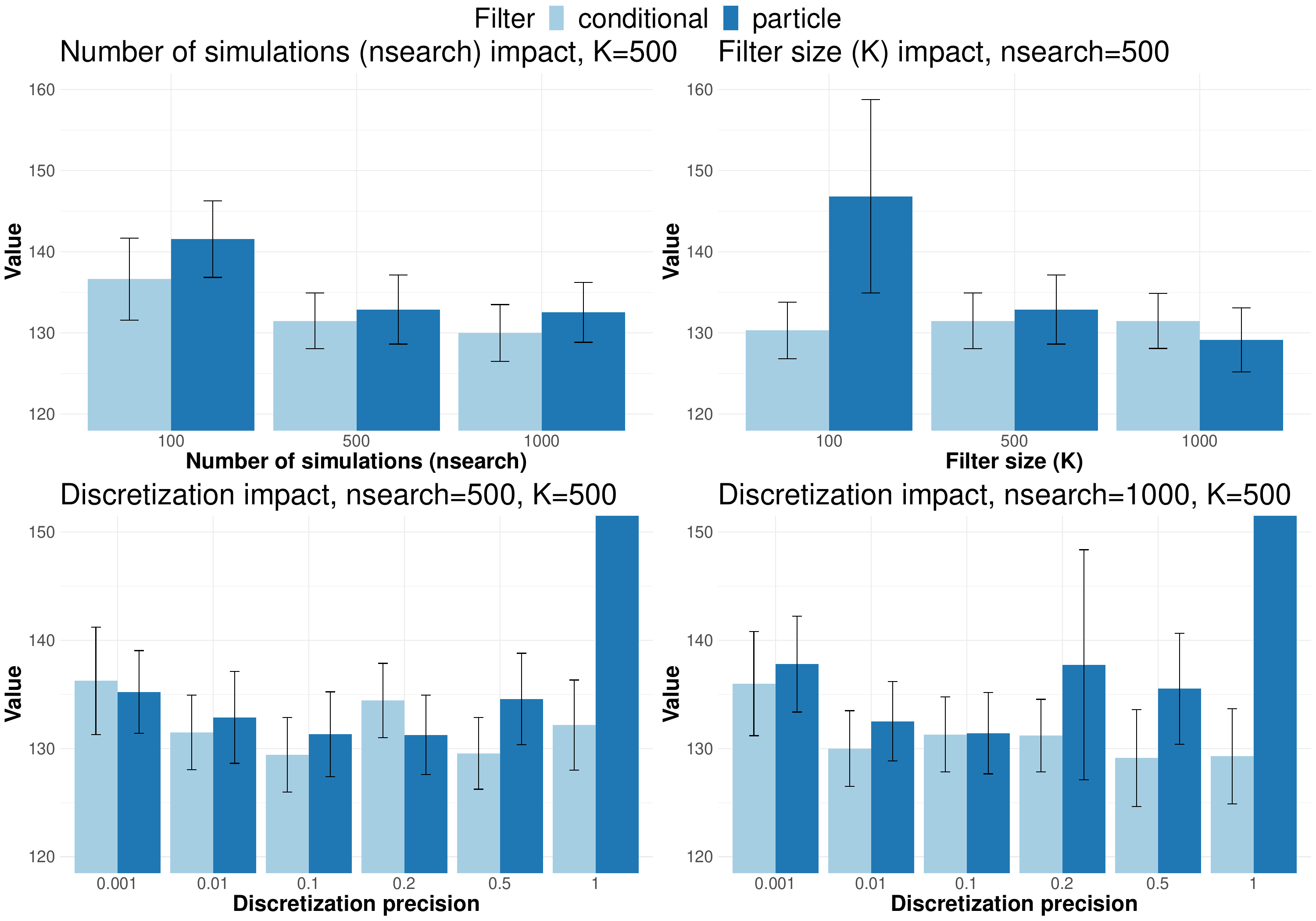}
\caption{{\bf Impact of POMCP parameters on the estimated value function.} Top left: increasing the number of simulations for filters with $500$ initial states improves the average trajectory costs. Top right: increasing the number of atoms in the filter improves the performance of the particle filter but not the conditional filter.} \label{fig:simus}
\end{figure}

\paragraph{Increasing the number of exploratory simulations yields the best performance gain.}
 For a fixed number of sampled belief states $K$ and precision $\mathcal{D}$, increasing the number of simulations $n_{\text{search}}$ improved the performance of the algorithm while decreasing the variance in the simulations for both filter types (see the top left panel of Figure~\ref{fig:simus} for $K=500$, $\mathcal{D}=0.01$). 
In the case of the conditional filter, the runtime increases linearly, from $10^3$ seconds per trajectory with $n_{\text{search}}=100$ simulations to $10^4$ seconds per trajectory with $n_{\text{search}}=1000$. In the case of the particle filter, the runtime is $3$ times as much for $100$ simulations since when $n_{\text{search}}<K$ additional simulations have to be performed to compute the particle filter at time $n+1$. The difference decreases to only $1.2$ times the runtime of the conditional filter for $n_{\text{search}}=1000.$ %

\paragraph{The particle filter requires a large belief state to achieve high performance.} As expected, for a fixed number of exploratory simulations $n_{\text{search}}$ and precision $\mathcal{D}$, increasing the number $K$ of particles in the particle filter led to a tremendous improvement for the particle filter together with a significantly decreased variance in the trajectory costs, while it had no impact on the conditional filter (see top right panel of Figure~\ref{fig:simus} for $n_{\text{search}}=500$, $\mathcal{D}=0.01$). Similarly, $K$ has no runtime impact for the  conditional filter, while it leads to exponential increase of the runtime for the particle filter.

\paragraph{POMCP for controlled PDMPs requires one additionnal tuning: the precision of the tree observation nodes.} The bottom two panels of Figure~\ref{fig:simus} illustrate the impact of the precision $\mathcal{D}$ for a fixed filter size ($K=500$) and two different numbers of simulations ($n_{\text{search}}=500$, left, and $n_{\text{search}}=1000$, right). For a smaller number of simulations, decreasing the precision improves the result up to $\mathcal{D}=0.1$, and then worsens them again. This tendency is still observed for the particle filter when $n_{\text{search}}$ increases, while the performance of the conditional filter is optimal with the loosest precision, $\mathcal{D}=1$. 
Those results are the consequence of two factors: as detailed in the Methods Section, each simulation creates a novel node in the tree exploration, where a node is a set of potential future trajectories with their estimated costs. When the precision is very fine, each simulation produces a different future observation and hence the estimation relies on one-step forward simulations which may miss future events. On the other hand, when the precision is very loose, each simulation will build on the previous simulation to explore a step forward, yielding very long trees with few branches. This will also miss the variability of different outcomes.
The second factor comes from the way the filters are constructed. Particle filters rely on comparing simulations to observations. When the precision is very loose, almost all simulations will be accepted, creating a strong bias in the belief of the current state, that will propagate from time-point to time-point. 
As the conditional filter update does not rely on simulations, hence neither on precision, there is no propagation of uncertainty from step to step, and when the number of exploratory simulations $n_{\text{search}}$ is large enough to guarantee some diversity in the tree exploration, estimating the cost of each decision from longer trajectories will provide better results. \\

Finally, one may note that the gain in using conditional filters is mostly apparent in extreme parameter scenarios, for instance with very low number of particles $K$, with very high precision rates, etc. Provided the user has enough computing budget, both filters tend to provide very similar results.

\subsubsection{Study 2: Adapted POMCP outperforms  the dynamic programming approach}

In this study we compare the results of three resolution strategies calibrated with their optimal parameters on biological relevant outcomes: the death rate, the \gls{pfs}
time, that is the time from entry in the study to the first relapse, the time spent under treatment, the number of visits to the hospital, and the cost. 
Those quantities were normalized so that they range between 0 and 1, and such that an optimal result is $0$. To do so, death rate was normalized so that a random treatment strategy yields 1 (here 5\% of patients); the PFS was transformed as 1-(PFS/H) (where we recall that $H$ is the study horizon), so that a patient who does not relapse has normalized PFS equal to $0$; the time spent under treatment was normalized by $H$, the number of visits was normalized as $\frac{N_{\text{Visit}}-40}{160-40}$ since over the horizon, a visit every $15$ days produces $160$ visits, whereas a visit every $60$ days produces $40$ visits; and finally the cost was normalized as $\frac{C-v_0}{C_{\text{random}}-v_0}$ where $v_0$ is the the best approximation of the optimal value obtained through discretizations in \cite{CdS23}, and $C_{\text{random}}$ is the average cost of the random strategy. 

\begin{figure}[ht!]
\centering
\includegraphics[width=8cm]{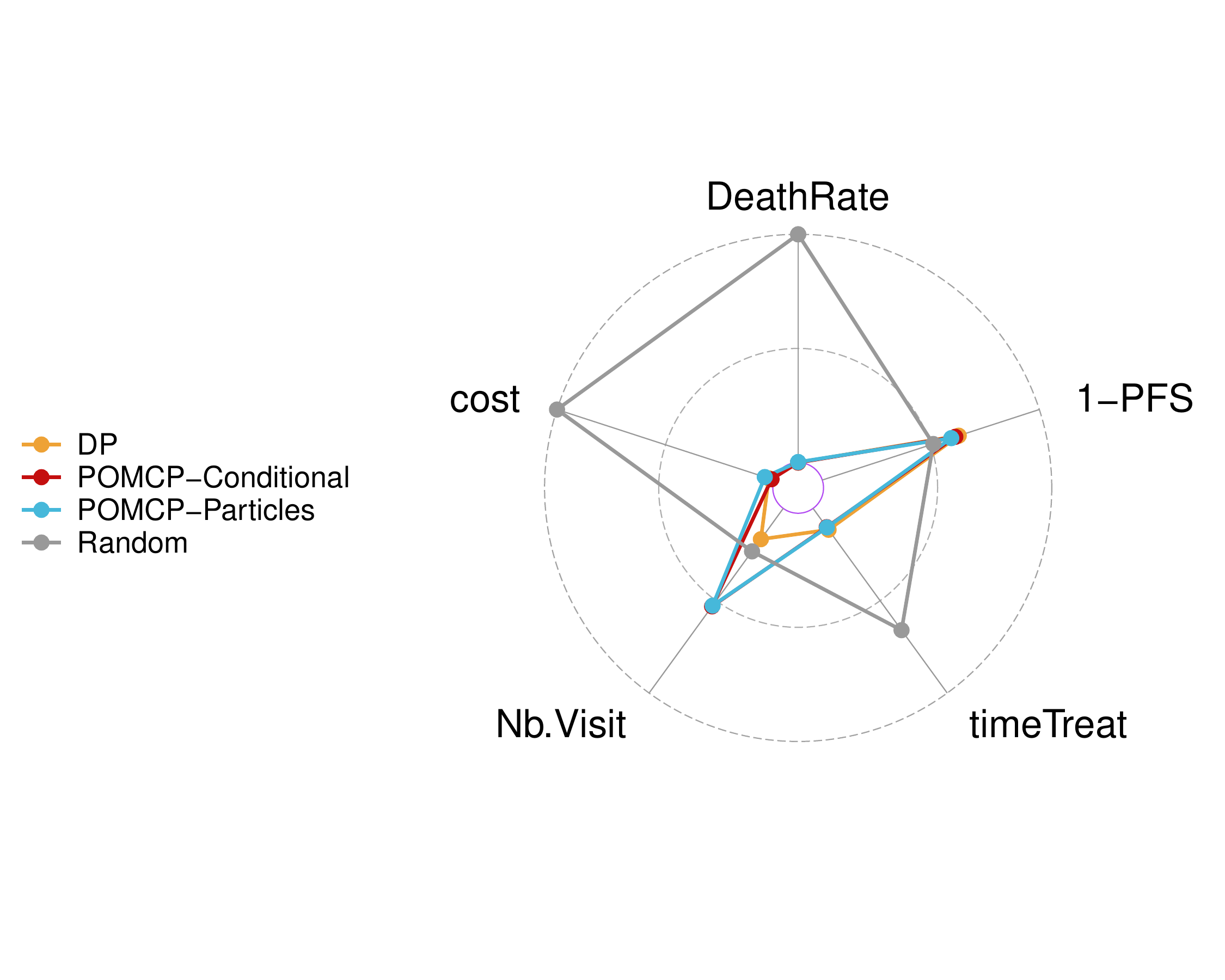}
\includegraphics[width=5cm]{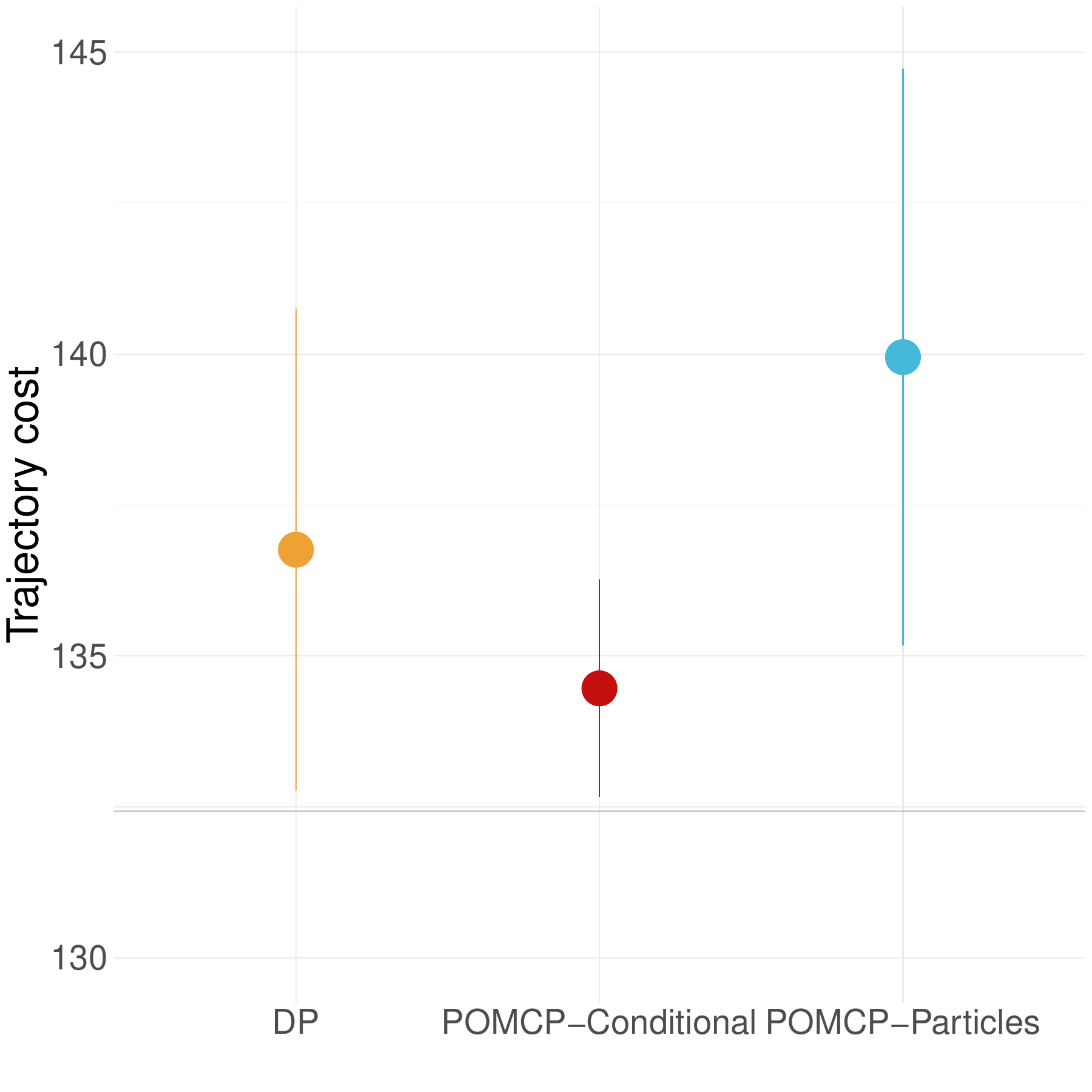}
\caption{\textbf{a)} Radar plot comparing performances of 4 solution strategies on death rate, progression-free survival (PFS), time spent under treatment, average number of visits per patient, and average trajectory cost. An optimal strategy would be the inner-circle. \textbf{b)} Barplot of trajectory cost for 500 simulations under three main strategies: POMCP with particle filter, POMCP with conditional filter and Dynamic Programming on discretized processes.} \label{fig:radar}
\end{figure}

Here again, we simulated $500$ trajectories with each strategy under the same cost parameters. The results are summarized in the Radar plot of Figure \ref{fig:radar}, and additional visual information on average trajectory cost are given in the barplots. In the Radar plot representation, a perfect strategy should delimitate the inner circle. 

Compared strategies are the discretization/DP approach (DP) from \cite{CdS23} that relies on exact resolution by dynamic programming of the discretized POMPD, the adapted POMCP with the conditional filter (POMCP-Conditional) and the adapted POMCP with the particle filter (POMCP-Particles).
Interestingly, the combination of POMCP with the conditional filter yields the lowest average trajectory cost and the shortest average time spent under treatment, by slightly increasing the number of visits and reducing PFS compared to the DP approach. However, the particle-based POMCP approach (relying fully on simulations with no other exploitation of the underlying model) yields cost almost as good as the previous two approaches, with increased number of visits but longest progression-free survival time. Importantly, out of $500$ simulations, one of the trajectories ended with a patient dying.

\section{Discussion}

In this paper, a mechanistic PDMP model for cancer evolution and treatment has been presented. PDMPs are very flexible tools that allow to model routine screening data with very few parameters with biological meaning. When embedded in a control framework, PDMPs usually suffer from the need to compute intractable integrals and thus resort to several layers of approximations with heavy computational burden. In our case, 
there are two major difficulties in solving the optimization problem for the controlled PDMP. The first one is related to the partial observations of the process, since the practitioner only observes some noisy measurement of the marker at visit dates, and not the overall condition of the patient nor the relapse dates. The second one comes from the continuous state space and continuous time dynamics of the process, which prevent direct use of exhaustive exploration solution strategies such as dynamic programming \cite{BT96}. 

In a previous work \cite{CdS23}, we have dealt with those difficulties by defining an equivalent fully observable Markov decision process on an enlarged state space through the use of \textit{conditional filters}, the conditional distributions of the hidden process given the observations. 
Then we discretized the state space of the original process, in order to obtain \textit{finite support} filters and discretized again these finite support approximate filters to obtain finitely many (belief) states. The fully discretized model can then be solved by dynamic programming.

Here we investigated another original solution approach exploiting filter objects under a different (simulation-based) dimension reduction strategy. We show that the inherent generator function of the PDMP can be exploited to make use of simulation-based solution strategies such as POMCP with excellent performance. Provided the number of simulations is large enough (either to explore the outcome space, or to construct consistent belief states), this approach can even outperform discretization approaches that exploit the knowledge of the underlying model. 
We have also proposed to combine both approaches, using discretization based conditional filters and simulation based solution stategy, resulting in a more robust algorithm (in particular less sensitive to the choice of POMCP parameters  and with more stable variance), but with little performance gain. 

The main advantage of the discretization/DP approach is that solutions are pre-computed for all new patients. This is especially useful under the assumption that all patients have the same dynamics with the same parameters. Its main drawback is that the model presented here is at maximum complexity for such an approach. In particular, it will become intractable if one wants to take into account more disease markers or more modes and treatments.

Conversely, the simulation-based approach can extend to any complexity provided simulations can be performed easily and fast. In addition, one of the main advantages of 
simulation-based approaches such as that presented here is that cost parameters can be modified with each new patient, 
yielding a patient-based procedure closer to precision medicine. This remains quite theoretical, as in practice calibrating cost parameters is a very difficult task, but with experience practitioners may be able to encode personal preferences, such as shorter life with better quality, or longer life at the price of more treatments, etc. 

In this work, we have adapted the POMCP algorithm to solving a controlled PDMP with {\em known} model, too complex to be solved through exact dynamic programming.
We made the assumption that the patient-disease model was known, which is a daring assumption. 
Therefore, one next step of our approach is to extend it by considering an unknown model and applying {\em Reinforcement Learning} methods \cite{SB98}.
A fundamental problem in Reinforcement Learning is the difficulty of deciding whether to select actions in order to learn a better model of the environment, or to exploit current knowledge about the rewards and effects of actions \cite{Katt17}. 
This is especially true in disease control problems.

\section{Methods}
\subsection{Datasets, parameters and code availability}
In order to propose a simulation study as realistic as possible we have used real data to infer the parameters of the design. The data come from the follow-up of 748 multiple myeloma patients registered in the $2009$ IFM clinical trial described in the Results Section. An example of data is given in Figure \ref{fig:patdata}. From this data, we opted for the exponential form of the dynamics in the disease states with boundaries $\zeta_0=1$ and $D=40$ for remission and death levels, simply calibrated as the minimal and maximal values in the data set. Then, as described in the Results Section, $3$ parameters had to be calibrated, and all hyperparameters are explicitly given in the github repository \cite{pomcp4pomdp}.

\begin{figure}[ht!]
\centering
\includegraphics[width=12cm]{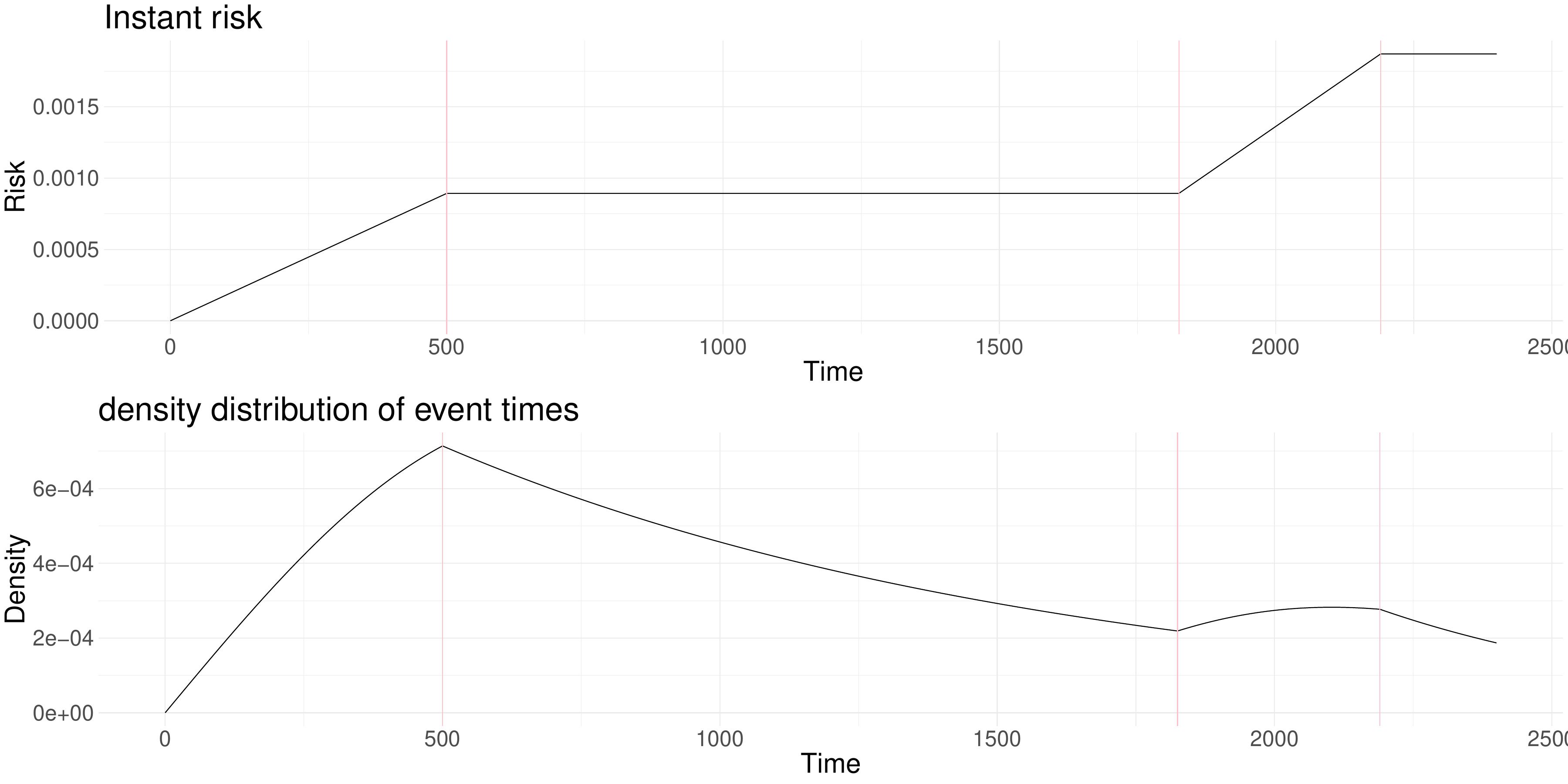}
\caption{Risk and Density functions for standard relapse from remission condition to disease $b$ condition (similar shapes for standard relapse to disease $a$).} \label{fig:risk}
\end{figure}

For the risk function $\lambda$, we choose to distinguish the \textit{standard relapse} (from remission to disease state) from the \textit{therapeutic escape} (from a disease state under appropriate treatment to the other disease state). We then further separate the risk by disease and treatment, so that for any treatment $\ell$ and any state $s=(m,\zeta,u)$, one has $\lambda^\ell(s)=\lambda_m^\ell(\zeta,u)$, where the form of $\lambda_m^\ell$ is specified in Table \ref{tab:intensite}. 
\begin{table}[tp]
\caption{Risk function of the controlled continuous time PDMP.}
\label{tab:intensite}
\centering
\scalebox{0.8}{
\begin{tabular}{l|lll}
&$\ell=\emptyset$&$\ell=a$&$\ell=b$\\
\hline
$m=0$&$\lambda_m^\ell(\zeta,u)=(\mu_1\!+\!\mu_2)(u)$&$\lambda_m^\ell(\zeta,u)=\mu_2(u)$&$\lambda_m^\ell(\zeta,u)=\mu_1(u)$\\
$m=1$&$\lambda_m^\ell(\zeta,u)=0$&$\lambda_m^\ell(\zeta,u)=\mu'(\zeta)$&$\lambda_m^\ell(\zeta,u)=0$\\
$m=2$&$\lambda_m^\ell(\zeta,u)=0$&$\lambda_m^\ell(\zeta,u)=0$&$\lambda_m^\ell(\zeta,u)=\mu'(\zeta)$ \\
$m=3$&$\lambda_m^\ell(D,u)=0$&$\lambda_m^\ell(D,u)=0$&$\lambda_m^\ell(D,u)=0$
\end{tabular}}
\end{table}
 For the \textit{standard relapse}, the risk $\mu_i$ for disease $i$ was chosen as piecewise increasing linear functions calibrated such that the risk of relapsing increases until some duration $\tau_1$ (average of standard relapses occurrences), then remains constant, and further increases between say $\tau_2$ and $\tau_3$ years (to model late or non-relapsing patients). This function and corresponding density are illustrated in Figure \ref{fig:risk} for disease $b$.

For the therapeutic escape, we chose to fit a Weibull survival distribution of the form
\begin{align*}
\mu'(\zeta)&=(\tilde \beta\zeta)^{\tilde \alpha},
\end{align*}
with $-1<\tilde\alpha<0$ to account for a higher relapse risk when the marker decreases. We arbitrarily chose $\tilde\alpha=-0.8$ and calibrated $\tilde\beta=1000$ such that only about 5\% of patients experience a therapeutic escape. 

The aggressiveness of the disease/treatment efficiency $v/v'$ may depend on both treatment $\ell$ and mode $m$. Specific values are thus denoted by $v_m^\ell/ v_m'$, see Table \ref{tab:flot}.
\begin{table}[tp]
\caption{Flow $\Phi$ of the controlled continuous time PDMP. For any treatment $\ell$, initial state $s=(m,\zeta,u)$, and duration $t$, $\Phi^\ell(m,\zeta,u,t)=(m,\Phi^\ell_m(\zeta,t),u+t)$ is the state of the patient after a time $t$ starting from state $s$ at time $0$ if no change of condition or treatment occurred.}
\label{tab:flot}
\centering
\scalebox{0.8}{
\begin{tabular}{l|lll}
&$\ell=\emptyset$&$\ell=a$&$\ell=b$\\
\hline
$m=0$&$\Phi_m^\ell(\zeta,t)=\zeta$&$\Phi_m^\ell(\zeta,t)=\zeta$&$\Phi_m^\ell(\zeta,t)=\zeta$\\
$m=1$&$\Phi_m^\ell(\zeta,t)=\zeta e^{v_1^\emptyset t}$&$\Phi_m^\ell(\zeta,t)=\zeta e^{-v'_1 t}=\zeta e^{v^a_1 t}$&$\Phi_m^\ell(\zeta,t)=\zeta e^{v_1^b t}$\\
$m=2$&$\Phi_m^\ell(\zeta,t)=\zeta e^{v_2^\emptyset t}$&$\Phi_m^\ell(\zeta,t)=\zeta e^{v_2^a t}$& $\Phi_m^\ell(\zeta,t)=\zeta e^{-v'_2 t}=\zeta e^{v^b_2 t}$\\
$m=3$&$\Phi_m^\ell(\zeta,t)=\zeta$&$\Phi_m^\ell(\zeta,t)=\zeta$&$\Phi_m^\ell(\zeta,t)=\zeta$
\end{tabular}}
\end{table}
After setting aside patients that do not relapse (about 20\%), we first estimated remission and relapse times using maximal slope difference, and then fitted exponential regression models on each segment. We then clustered patients based on their relapse coefficient and chose the number of clusters using the slope heuristic on residual sum of squares. We obtained two groups and used the average values to obtain $v_1^\emptyset=0.02$ ($22$\% of patients) and $v_2^\emptyset=0.006$ ($78$\% of relapsing patients). We then computed the average of the treatment parameters for each group and obtained $v'_1=0.077$ and  $v'_2=0.025$. The data do not present patient relapsing under treatment, and therefore we could not estimate $v$ for therapeutic escapes or inappropriate treatments. Because we assume the aggressiveness in those circumstances should be smaller than under standard relapse,  we chose $v_1^b=0.01$ and $v_2^a=0.003$.

By separating the risk $\lambda$ by disease, the kernel function $Q$ is automatically fitted, since we assume that the marker level does not jump at relapses, and the mode is selected by the risk clocks leading to the jump, whichever rings first, see Table \ref{tab:Q}.
\begin{table*}[tp]
\caption{\label{tab:Q} Markov transition kernel $Q$ for the controlled PDMP. matrix $Q_m^\ell(m')$ for the controlled PDMP. For any treatment $\ell$ and initial state $s=(m,\zeta,u)$ a jump sends the patient to state $s'=(m',\zeta',u')$ sampled from the distribution $Q(\cdot|s,\ell)$. The following constraints are satisfied: $\zeta'=\zeta$, $u'=0$ and $m'$ is sampled from the discrete distribution $Q_m^\ell$.}
\centering
{
\begin{tabular}{l|l}
&$\ell=\emptyset$\\
\hline
$m=0$&$Q^\ell_m(m')=\1_{(m'\in\{1,2\})}\frac{\mu_{m'}(u)}{\mu_1(u)+\mu_2(u)}$\\
$m=1$&$Q^\ell_m(m')=\1_{(m'=0)}$ (possible only if $\zeta=\zeta_0$)\\
	   &$Q^\ell_m(m')=\1_{(m'=3)}$ (possible only if ${\zeta=D}$)\\
$m=2$&$Q^\ell_m(m')=\1_{(m'=0)}$ (possible only if $\zeta=\zeta_0$)\\
	   &$Q^\ell_m(m')=\1_{(m'=3)}$ (possible only if ${\zeta=D}$)\\
\hline
\hline
&$\ell=a$\\
\hline
$m=0$&$Q^\ell_m(m')=\1_{(m'=2)}$\\
$m=1$ &$Q^\ell_m(m')=\1_{(m'=2)}$ (possible only if $\zeta>\zeta_0$)\\
	  &$Q^\ell_m(m')=\1_{(m'=0)}$ (possible only if $\zeta=\zeta_0$)\\
$m=2$
	   &$Q^\ell_m(m')=\1_{(m'=3)}$ (possible only if ${\zeta=D}$)\\
\hline
\hline
&$\ell=b$\\
\hline
$m=0$&$Q^\ell_m(m')=\1_{(m'=1)}$\\
$m=1$&$Q^\ell_m(m')=\1_{(m'=3)}$ (possible only if ${\zeta=D}$)\\
$m=2$ &$Q^\ell_m(m')=\1_{(m'=1)}$ (possible only if $\zeta>\zeta_0$)\\
	  &$Q^\ell_m(m')=\1_{(m'=0)}$ (possible only if $\zeta=\zeta_0$)\\
\end{tabular}}
\end{table*}
Finally, we arbitrarily selected a centered Gaussian distribution with noise parameter $\sigma^2=1$ for the observation process. 

We resorted to extensive simulations study to select cost parameters that seemed reasonable (very few patients dying over our study horizon, on average not more that a fourth of the follow-up time spent under treatment, etc). We arbitrarily fixed the visit cost $C_V$ to $1$. We then fixed the death cost $M$ to $110$ so that with visits every $15$ days and early relapse, a patient would rather die that spend the entire horizon under treatment with numerous visits. 
We then fixed $\beta=0.1$ so that the penalty of applying an unnecessary treatment would count as $1.5, 3$ or $6$ times the visit cost depending on the choice of next visit date $r$. Finally, we selected $\kappa=1/6$ from extensive simulations so that for low marker observations (typically when it is hard to decide between relapse or remission with high noise) it might be preferable to wait for new data acquisition rather than treat by default. All codes and parameters are available at \url{https://github.com/acleynen/pomcp4pdmp} \cite{pomcp4pomdp}.

\subsection{Reminder on controlled PDMPs}
Here is a description of how to simulate a trajectory of a controlled PDMP between two consecutive visits to the medical center.
{The controlled PDMP parameters are $\lambda$, the disease risk function (distribution of duration until the next jump i.e. condition change), $Q$, the Markov kernel, defining the stochastic transition to the state reached after the next jump and $\{v^{\ell}_m\}$\footnote{For the sake of unified notation, we denote here $v_1^a=v'_1$, $v_2^b=v'_2$, and $v_0^\emptyset=v_0^a=v_0^b=0$}, the parameters of the exponential deterministic behavior of the marker between two jumps, defined from the current mode and treatment applied.}

\begin{algorithm}
\caption{Simulation of a trajectory between two consecutive decision times of a controlled PDMP.\label{SimuPDMP}}
\begin{algorithmic}[1]

\Procedure{SimulatePDMP}{$m$, $\zeta$, $u$, $\ell$, $r$}
    \State $t\leftarrow 0$
    \State $v\leftarrow v_m^\ell$
    \While{$t<r$}
        \State $S\sim \lambda$
        \State $S \leftarrow \min\{S, t^*(m,\zeta,u,\ell)\}$
        \If{$t+S>r$}
            \State \Return{$m$, $\zeta\exp(vr), u+r$}
        \Else
            \State $t\leftarrow t+S$
            \State $\zeta \leftarrow \zeta\exp(vS)$
            \State $u\leftarrow u+S$
            \State $m\sim Q(\cdot | m,\zeta,u,\ell)$
            \State $u\leftarrow 0$
            \State $v\leftarrow v_m^\ell$
        \EndIf
    \EndWhile
\EndProcedure
\end{algorithmic}
\end{algorithm}

Procedure {\sc SimulatePDMP} takes as input an initial position $X_t=s=(m,\zeta,u)$ with mode $m$, marker level $\zeta$, time since the last jump $u$, and a decision $d=(\ell,r)$ with treatment to be applied $\ell$ for a duration $r$ until the next visit to the medical center and returns the state $X_{t+r}=s'=(m',\zeta',u')$ of the process at time $t+r$ given that treatment $\ell$ was applied. 
At line $5$, $S\sim \lambda$ means that $S$ is sampled from the distribution with risk function $\lambda$, which means that it has the following survival function
$$\mathbb{P}(S>t)=e^{-\int_0^t \lambda(m, \zeta\exp(v\tau), u+\tau)d\tau}.$$
At line $6$, $t^*(m,\zeta,u,\ell)$ is the (deterministic) time to reach either the nominal value $\zeta_0$ or the death level $D$ from the current point $(m,\zeta,u)$ (if no change of condition or treatment occurs). The third variable $u$ representing the time since the last jump 
allows transitions described in {\sc SimulatePDMP} to be Markovian, i.e. to be independent of the previous transitions.

\subsection{Reminder on Partially Observed Monte-Carlo Planning}
In this section we give a description of the {original} POMCP algorithm\footnote{A Python implementation of the POMCP algorithm can be found here: \url{https://github.com/GeorgePik/POMCP}}.
The algorithm is called iteratively at each observation step $n$ and requires several entries in order to output a decision $d_n$ to be used by the operator. 
The inputs include a set of possible decisions (denoted $A$), a history \gls{ghn}  of successive observations and decisions up to step $n$, including the current observation $\omega_n$, an approximate (particle) \textit{filter} $\Theta^p(h_n)$ with support $B^p(h_n)$
, a simulator $\mathcal{G}(s,d)$ of state and observation at step $n+1$ together with their cost given a state $s$ and decision $d$ at step $n$, a stopping criterion {\sc Timeout}
and an arbitrary {\sc Rollout} strategy to provide a heuristic evaluation of an history, whenever needed.

\begin{algorithm}
\caption{Original POMCP algorithm \cite{SV10}}
\label{algo:POMCP}
\begin{multicols}{2}
\begin{algorithmic}[1]

\Procedure{POMCP}{$h_n$}
    \Repeat
        \If{$h_n=\emptyset$}
            \State{$s\sim \Theta^p_0$}
        \Else
            \State $s \sim \Theta^p(h_n)$
        \EndIf
        \State {\sc Simulate}$(s,h_n)$
    \Until{{\sc Timeout}()}
    \State $d^* \leftarrow \arg\min_d V(h_nd)$
    \State $V(h_n) \leftarrow \min_d V(h_nd)$
    \State \Return{$d^*$}
\EndProcedure

\Statex
\Procedure{Rollout}{$s$,$h$}
    \If{$s.t=H$}
    \State \Return{0}
    \EndIf
    \State $d\sim \pi_{\text{rollout}}(h)$
    \State $(s',\omega,c) \sim {\mathcal G}(s,d)$
    \State \Return{$c+${\sc Rollout}$(s',hd\omega)$}
\EndProcedure

\Statex
\Statex

\Procedure{Simulate}{$s$,$h$}
    \If{$s.t \geq H$}
        \Return{0}
        \EndIf
    \If{$h\not\in \mathcal{T}$}
        \ForAll{$d\in A$}
            \State $\mathcal{T}(hd) \leftarrow \langle N_{\text{init}}, V_{\text{init}}, \emptyset\rangle$
        \EndFor
        \State $C\leftarrow${\sc Rollout}$(s,h)$
        \State \Return{C}
    \EndIf
    \State $d^* \leftarrow \arg\min_d V(hd) - \alpha \sqrt{\frac{\log\left(N(h)\right)}{N(hd)}}$
    \State $(s',\omega,c) \sim {\mathcal G}(s,d^*)$
    \State $C\leftarrow c+${\sc Simulate}$(s',hd^*\omega)$
    \State $B^p(h)\leftarrow B^p(h)\cup\{s\}$
    \State $N(h)\leftarrow N(h)+1$
    \State $N(hd^*) \leftarrow N(hd^*)+1$
    \State $V(hd^*) \leftarrow V(hd^*)+\frac{C-V(hd^*)}{N(hd^*)}$
    \State \Return{C}
\EndProcedure
\end{algorithmic}
\end{multicols}
\end{algorithm}

The POMCP algorithm involves several data structures: 
\begin{itemize}
    \item Simulated states $s =(m, \zeta, u)$. 
    \item Decisions $d = (\ell, r)$, belonging to a finite decision space $A$, preferably small. 
    \item Observations $\omega = ({y}, t)$. The original algorithm assumes that they belong to a finite observation space $\Omega$ of limited size. 
    \item Histories $h = \langle \omega_0 d_0 \omega_1 d_1 \cdots, d_{n-1} \omega_{n} d' \omega' d" \omega" \dots \rangle$.
    Histories represent sequences of decisions and observations of variable lengths.
    They are the concatenation of the sequence of past observations/decisions $h_n = \langle \omega_0 d_0 \omega_1 d_1 \cdots d_{n-1} \omega_{n}\rangle$ plus an arbitrary sequence of {\em future} observations/decisions, built using the {\em rollout strategy} and a simulation model of the POMCP  (line 28, $h\leftarrow hd^*\omega$).
    \item $\mathcal T$ is a tree data structure rooted at the initial history $h_n$. Each node of $\mathcal T$ will correspond to an extended history, ended by either a decision or an observation.  
    Each simulation step creates a novel node in the tree, and histories $h$ are attached to $\mathcal T$ by appending the corresponding novel decision/observation to the parent node's history. In addition, we also attach to each node (now denoted ${\mathcal T}(h)$ or ${\mathcal T}(hd)$) 
    (i) integer numbers $N(h)$ or  $N(hd)$ where $N(h)$ corresponds to the number of times  the history $h$ has been simulated and  $N(hd)$ to the number of times $d$ has been selected after $h$ was encountered, 
    (ii) real numbers $V(hd)$ corresponding to an estimate of the {\em cost value} of $hd$, that is the expected sum of future costs
    obtained if the optimal policy is applied after $h$ has been encountered and $d$ selected, until the final decision step  and 
    (iii) $\Theta^p(h)$, called a particle filter, which is a discrete uniform distribution on a set $B^p(h)$ of states $s$ compatible with the current history $h$. 
\end{itemize}

When {\sc Simulate} is called with entry a history $h$ that already belongs to $\mathcal T$, it updates the values of $N(h)$ and $B^p(h)$ as well as the values $N(hd)$ and $V(hd)$ for all the successor nodes\footnote{It is a property of the algorithm that whenever $h\in {\mathcal T}$, $hd\in\mathcal T$ as well.} $hd$.
When {\sc Simulate} is called with entry a history $h$ which does not yet belong to $\mathcal T$ (as is the case initially for $h_n$), it appends $h$ as well as all its successor nodes $hd$ to $\mathcal T$ and initializes their values $N(h), B^p(h), N(hd)$ and $V(hd)$. 

Procedure {\sc Simulate} is based on a {\em generator function}, $(s', \omega, c)\sim {\mathcal G}(s, d)$ that generates a successor (hidden) state $s'$, an observation $\omega$ and an immediate cost $c$, from decision $d$ applied in current (hidden) state $s$.
Repeated calls to $\mathcal G$ are used to progressively expand  $\mathcal T$.

Simulation sequences and updates are performed, starting from $h_n$, until {\sc Timeout}() function requires to stop (generally, after an arbitrary number $n_{\text{search}}$ of trajectories have been simulated or a fixed amount of time has been spent). 
Then, the decision $d^*$ which maximizes $V(h_nd^*)$ is applied to the real-world system, and a real-life observation $\omega\in\Omega$ is obtained.
The new real-world history becomes $h_{n+1} = h_nd^*\omega$ and $\mathcal T$ is pruned, so that the new tree is rooted\footnote{The interest of pruning $\mathcal T$ instead of starting with an empty tree in $h_{n+1}$ is to exploit past simulations in the computation of the next decision.} in $h_{n+1}$. \newline

POMCP proposes strategies to select the input rollout and filters when the user has no knowledge on the process. A typical rollout strategy may simply involve selecting the decisions randomly from the set $A$. 
The following particle filter update procedure, included in procedure {\sc Simulate}, is suggested: 

\begin{itemize}
    \item If $h=\emptyset$, sample $s\sim \Theta^p_0$. In the initial step of the algorithm, we simulate random particles from an arbitrary belief state. 
    \item If $h\neq\emptyset$, sample $s\sim\Theta^p(h)$ where $\Theta^p(h)$ is the uniform discrete distribution on the finite nonempty set $B^p(h)$. Indeed, if $h\neq\emptyset$, this means that procedure {\sc Simulate}$(s',h)$ has already been called at least once for some $s'$ and thus $B^p(h)\neq\emptyset$ ($B^p(h)$ contains at least $s'$).

    \item Sample $(s',\omega,c) \sim {\mathcal G}(s,d^*)$. This step is performed in line 27 of the POMCP algorithm.
    \item if $|\omega-\omega'|=0$, update $B^p(hd^*\omega)\leftarrow B^p(hd^*\omega)\cup\{s'\}$.
\end{itemize}
As the number of samples increases, the supports of the filters $B^p(h)$ will contain more and more particles and converge to the empirical distribution $\mathbb P(\cdot|h)$ of hidden states given the observed trajectory $h$. When the state space is finite, $\Theta^p(h)$ can be seen as an histogram approximation of $P(\cdot|h)$.

\subsection{Adapted POMCP algorithm to the case of controlled PDMPs}
It seems natural to apply a POMCP algorithm to optimize a PDMP control strategy in the context of disease control. Indeed, PDMPs have natural {\em generator} functions since algorithm {\sc SimulatePDMP} can be naturally augmented with an observation simulator and a cost function, in order to obtain a generator function $\mathcal{G}$ as described above. 
We propose three adaptations of the original POMCP algorithm to exploit the particular framework of controlled PDMPs.

\paragraph{Rollout policies} The original POMCP algorithm describes the possible rollout policies as  {\em admissible} policies, meaning that actions choices should only depend on the history  of past actions and observations. Instead, we exploit here the PDMP generator which provides both hidden states and observations. This allows to design rollout policies exploiting hidden states instead of noisy observations. In our medical framework, 
one may build interesting {\em rollout policies} exploiting the hidden mode of the disease to provide good heuristics to the simulation part of POMCP.
In practice, we propose to compare the two following rollout policies:
\begin{enumerate}
    \item The (admissible) uniform  policy : $\pi_{\text{unif}}(\omega=(y,t)) \sim \mathcal{U}(\{\emptyset, a, b\}\times\{15, 30, 60\})$  
    \item The (non-admissible) mode policy : $\pi_{\text{mode}}(s=(m,\zeta,u))=\pi_{mode}(m)= \begin{cases}
    \{\emptyset,15\}\; \text{ if }\; m=0,\\
    \{a,15\}\; \text{ if }\; m=1,\\
    \{b,15\}\; \text{ if }\; m=2.    
    \end{cases}$
\end{enumerate}

The mode policy, being based on the full observation of the process, is likely to underestimate the real cost of an optimal control policy, which is a useful property for the convergence of a heuristic search method \cite{RN02}. 

In our simulation study we observed that a mode-based rollout policy can be particularly efficient.  

\paragraph{Observation space} The time, state and observation spaces of the disease control PDMP model are continuous. This means that the probability of simulating exactly the same history $h$ twice is zero if we apply the POMCP procedure as such. Thus, the tree depth and the size of filters supports $B^p(h)$ may never exceed 1. This is particularly annoying since the POMCP algorithm convergence proof only holds when $\Theta(h)$ is close to the true empirical belief state, which (approximately) holds when the size of the support $B^p(h)$ tends to $+\infty$. Indeed, the {\sc POMCP} procedure in \cite{SV10} requires that $B^p(h)$ contains at least $K$ particles, when $h$ is non empty.
In practice, for the controlled PDMP case, observations are made of pairs $(y, t)$ of a continuous-value observed marker level and discrete time of current decision. 
Therefore, we discretize the observation space into a set of contiguous intervals and group together observations belonging to the same interval. 

The continuous nature of the model also prevents the exact computation of the filter, hence we resort to the use of particle or conditional filters. The construction of the former is slightly adapted from the initial {\sc POMCP} algorithm to fit our model. Indeed, as the true process still produces continuous-valued observations, the last action of the particle filter update procedure is modified to updating $B^p(hd^*\omega)$ with $s'$ only if $|\omega-\omega'|<{\mathcal D}$ (with ${\mathcal D}$ chosen by the user, typically of the same magnitude as the discretization precision). It may still happen that this procedure selects only states with wrong mode $m$ at step $n$ (i.e. $\Theta^p(h_n)$, which is only an approximation of the true filter, does not contain the true hidden mode $m$), and hence cannot generate compatible states at step $n+1$ due to diverging dynamics of the process in the different modes.
This is not a simple matter of statistical accuracy, but a very practical problem. 
When this happens, we say that $B^p(h_{n+1})$ is {\em deprived} of particles and we cannot go on applying {\sc POMCP} to $h_{n+1}$.

This problem of particle deprivation can be mitigated by a few modifications of the standard particle filter construction: 
\begin{itemize}
    \item Assume that $B^p(h_{n+1})$ is empty or too small, whatever the number of simulations of $s_n\sim B^p(h_n)$ followed by a call to ${\mathcal G}(s_n,d_n^*)$ we perform. Then, we may go back in the history and resimulate $s_{n-1}\sim B^p(h_{n-1})$ followed by two successive calls to $\mathcal G$: 
    \begin{itemize}
        \item $(s_n, \omega', c') \leftarrow {\mathcal G}(s_{n-1}, d^*_{n-1})$ and, provided that $|\omega'- \omega_n|<{\mathcal D}$,
        \item $(s_{n+1}, \omega'', c'') \leftarrow {\mathcal G}(s_n,d^*_n)$, hoping that now, $|\omega''-\omega_{n+1}|<{\mathcal D}$.
    \end{itemize}
    A particle $s_{n+1}$ is then added to $B^p(h_{n+1})$ whenever the two above conditions are met.
    \item Assume that $B^p(h_{n+1})$ is non-empty but still too small ($|B^p(h_{n+1})|\ll K$) after the previous step was applied a large number of times. We may perform particle revigoration by resampling particles from $B^p(h_{n+1})$ and duplicate them.
    \item Finally, when everything fails, we may perform a large number of sampled transitions from $B^p(h_n)$ (e.g. $1000\times K$) and keep the $K$ particles $s'$ in the generated samples$(s', \omega', c')$ with minimal distance  $|\omega_{n+1}-\omega'|$, with some arbitrary distance definition. 
\end{itemize}

In our experimental studies, we applied these three modifications in turn, whenever needed, until we got belief states $B^p(h_n)$ of cardinality at least $K$.

\paragraph{Filters}
We propose to modify the original {\em particle filter} of POMCP to incorporate the {\em conditional filter}. Hence the algorithm is modified as follows: starting from an initial arbitrary belief filter $\Theta^c_0$ for $h_0=\emptyset$ as in the original POMCP, when the new history becomes $h_{n+1}=h_nd^\star\omega$, we compute the new filter $\Theta^c_{n+1}$ as a deterministic function of $\Theta^c_n$ and $d^\star,\omega$ (see \cite{CdS23} for its specific form) and sample $K$ particles from $\Theta^c_{n+1}$ to generate a set of plausible hidden states. 

Depending only on the current belief and the new observation (and not on simulations), this filter does not suffer from the propagation of approximations and particle deprivation. Moreover, the computational burden of the simulations is replaced by the computation of weighted sums which are particularly efficient in matrix programming languages. 

\section{Supporting information}
\subsection{State of the art}
\paragraph{Artificial Intelligence in medicine}
The use of artificial intelligence methods in medicine has recently exploded, as was shown for example in \cite{Kumar22}. Yet the vast majority of these works focus on diagnosis and prognosis, rather than treatment and follow-up.
There are a few studies related to treatment of diseases, though.
\cite{Bhinder21} provides a review of approaches to cancer treatment focused on medical decision support. 
In the context of {\em drug design} for cancer treatment, (deep) \gls{rl}
approaches have been proposed recently \cite{Olivecrona17,Popova18}, but with the disadvantage of being black-box approaches, preventing the practitioner from access to an explainable model of disease evolution and treatment. 
Uncertainty quantification in models of assisted decision making have been proposed to alleviate this problem \cite{Begoli19}.
In the same line of work, \cite{Benzekry20} advocate learning of mechanistic models of cancer evolution/treatment in order to help decision making.
However, the latter approach requires a large amount of data {\em prior} to applying any treatment action, in order to learn a model that may prove only partially valid as decisions influence the disease dynamics. 
The authors advocate that the approach may be rendered more efficient if learning phases are interleaved with actual decision phases.

\paragraph{An alternative view of controlled PDMPs}
Controlled PDMPs can be modelled as continuous space POMDPs, in the way we proposed in this paper. However, they can also be seen as a particular subclass of continuous-time (Partially Observed) {\em Semi-Markov Decision Processes} \cite{Howard64}. Fully-observed continuous-time Semi-Markov Decision Processes extend Markov Decision Processes by including random continuous durations of state transitions and by considering that decisions can only be made at transition times. Several reinforcement learning solution approaches have been proposed, both in the fully-observed \cite{BradtkeDuff95,Doya2000} and partially-observed \cite{Doshi2020} cases to solve these problems. An alternative approach to ours could be to cast the PDMP model of cancer treatment into the continuous-time SMDP model and look for specializations of the existing simulation based solution algorithms.

\subsection{Supplementary simulation results on POMCP parameters}
 Here we provide raw results for a series of parameter we tried to tune to optimize the POMCP algorithm, but that did not seem to bring additional improvement in our framework. Table \ref{table_resultats_1} shows the impact of the {trade off parameter $\alpha'$} in a few scenarios for the particle filter.
 
\begin{table}[ht!]
    \centering
        \caption{Raw results for the particle filter varying parameters $\alpha'$, $n_{\text{search}}$ and $K$. For each parameter set, $n=500$ trajectories were simulated. The Value column is the average cost of the trajectories over the $n$ trajectories, and $\hat{\sigma}$ its empirical variance. We also recorded the runtime of optimizing each trajectory (duration column).}
\label{table_resultats_1}
  \begin{tabular}{c|c|c|c|c|c|c|c|c}    \hline
     Filter & $\pi_{rollout}$ & $n_{\text{search}}$ & $K$ & {$\alpha'$} & Value & $1.96\hat{\sigma}/\sqrt{n}$ & duration & duration s.d  \\ \hline
    \hline
     particle & $\pi_{mode}$ & 100 & 100 & 0.2 & 161.93 & 14.79 &1730 & 982\\
   particle & $\pi_{mode}$ & 100 & 100 & 0.5 & 156.01 & 13.01& 1629 & 913\\
    particle & $\pi_{mode}$ & 100 & 100 & 0.8 & 165.63 & 13.51 & 1714 & 1053\\
    particle & $\pi_{mode}$ & 100 & 100 & 0.99 & 147.24 & 6.17 & 1641& 971\\
    \hline \hline
    particle & $\pi_{mode}$ & 100 & 500 & 0.2 & 141.10 & 9.83 &2727 & 683\\
   particle & $\pi_{mode}$ & 100 & 500 & 0.5 & 141.56 & 4.72& 2771 & 699\\
    particle & $\pi_{mode}$ & 100 & 500 & 0.8 & 134.98 & 4.44 & 2646 & 685\\
    particle & $\pi_{mode}$ & 100 & 500 & 0.99 & 133.57 & 3.56 & 2640& 611\\
    \hline \hline
    particle & $\pi_{mode}$ & 500 & 100  & 0.2 & 146.30&8.27 & 4617&660\\
    particle & $\pi_{mode}$ & 500 & 100 & 0.5 & 146.82&11.91&4047&670\\
    particle & $\pi_{mode}$ & 500 & 100 & 0.8 & 145.87&13.01&3708&610\\
    particle & $\pi_{mode}$ & 500  & 100 & 0.99& 140.91&6.35&3454&614\\
    \hline \hline
   particle & $\pi_{mode}$ & 500  & 500& 0.2 & {135.99}&4.00&5182&598\\
    particle & $\pi_{mode}$ & 500 & 500 & 0.5 & {132.88}&4.25&5068&643\\
    particle & $\pi_{mode}$ & 500  & 500 & 0.8 & {136.14}&8.08&5271&698\\
    particle & $\pi_{mode}$ & 500 & 500 & 0.99& {129.42}&5.06&4946&724\\
    \end{tabular}
\end{table}

 To allow adaptive selection of the trade off parameter $c$ we tried three dynamic procedures to exploit or explore more depending on our trust in the current patient state. To do so, we define the state entropy as $E_t=\sum_{m=0}^2 p_m \log(p_m)$, $p_j=\sum_{s=(m,\zeta,u)\in B} \mathbf{1}\{m=j\}$ and $E_{\max} =\log(1/3)$ and consider the following:
\begin{itemize}
    \item entropy: $\alpha_t=E_t/E_{\max}$
    \item rev-entropy: $\alpha_t=1-E_t/E_{\max}$
    \item rev-entropy-2: $\alpha_t=1-E_t/2E_{\max}$.
\end{itemize}
However, none of those procedures improved the results, as illustrated in selected examples in Table \ref{table_resultats_2}.

\begin{table}[h!]
    \centering
        \caption{Simulation results with adaptative choice of the exploration/exploitation parameter $\alpha'$. For each parameter set, $n=500$ trajectories were simulated. The Value column is the average cost of the trajectories over the $n$ trajectories, and $\hat{\sigma}$ its empirical variance. We also recorded the runtime of optimizing each trajectory (duration column).}
    \label{table_resultats_2}
    \begin{tabular}{c|c|c|c|c|c}
         Filter &  $n_{\text{search}}$ & $\alpha'$ & Value & $1.96\hat{\sigma}/\sqrt{n}$ & duration  \\ \hline \hline
        conditional & $100$ & entropy &$138.63$ & $5.06$ & $776$ \\
        conditional & $100$ & rev-entropy &$131.94$ & $3.70$ & $786$ \\
        conditional & $100$ & rev-entropy-2 &$133.75$ & $3.70$ & $770$ \\
\hline \hline        
        particles & $100$ & entropy & $142.17$ & $9.88$ & $2421$ \\
        particles & $100$ & rev-entropy & $143.27$ & $10.34$ & $2473$ \\
        particles & $100$ & rev-entropy-2 & $135.28$ & $3.82$ & $2438$ \\
\hline \hline
        conditional & $1000$ & entropy & $131.78$ & $4.70$ & $8313$ \\
        conditional & $1000$ & rev-entropy & $131.41$ & $3.45$ & $8432$ \\
        conditional & $1000$ & rev-entropy-2 & $132.73$ & $3.56$ & $8332$ \\
\hline \hline
        particles & $1000$ & entropy & $133.39$ & $3.74$ & $9994$ \\
        particles & $1000$ & rev-entropy & $135.64$ & $3.74$ & $10047$ \\
        particles & $1000$ & rev-entropy-2 & $131.89$ & $4.24$ & $10028$ \\
    \end{tabular}
\end{table}

\printglossaries  

\bibliographystyle{plain} 
\bibliography{biblio}

\end{document}